\documentclass{amsart}
\usepackage{amssymb}

\numberwithin{equation}{section}

\newtheorem{thm}{Theorem}[section]
\newtheorem{prop}[thm]{Proposition}
\newtheorem{lem}[thm]{Lemma}
\newtheorem{cor}[thm]{Corollary}
\theoremstyle{definition}
\newtheorem{defn}[thm]{Definition}
\theoremstyle{remark}
\newtheorem{rem}[thm]{Remark}

\newcommand{\Z}{\mathbb{Z}}
\newcommand{\Q}{\mathbb{Q}}

\newcommand*{\Satakepar}{\mathcal{S}}
\newcommand{\red}{\mathrm{red}}
\newcommand{\GL}{\mathop\mathrm{GL}\nolimits}

\newcommand{\Ord}{\mathrm{Ord}}
\newcommand{\Sp}{\mathrm{Sp}}
\newcommand{\trivrep}{\mathbf{1}}
\newcommand{\id}{\mathrm{id}}
\DeclareMathOperator{\cInd}{c-Ind}
\DeclareMathOperator{\Ind}{Ind}
\DeclareMathOperator{\supp}{supp}
\DeclareMathOperator{\Hom}{Hom}
\DeclareMathOperator{\End}{End}
\DeclareMathOperator{\Stab}{Stab}
\DeclareMathOperator{\Ker}{Ker}

\DeclareMathOperator{\Imm}{Im}
\DeclareMathOperator{\Gal}{Gal}

\title[On a classification of modulo $p$ representations]{On a classification of irreducible admissible modulo $p$ representations of a $p$-adic split reductive group}
\author{Noriyuki Abe}
\address{Creative Research Institution, Hokkaido University, N21 W10, Kita-ku Sapporo 001-0021, Japan}
\email{abenori@math.sci.hokudai.ac.jp}

\begin{document}
\begin{abstract}
We give a classification of irreducible admissible modulo $p$ representations of a split $p$-adic reductive group in terms of supersingular representations.
This is a generalization of a theorem of Herzig.
\end{abstract}
\maketitle

\section{Introduction}\label{sec:Introduction}
Let $p$ be a prime number and $F$ a finite extension of $\Q_p$.
In this paper, we consider modulo $p$ representations of (the group of $F$-valued points of) a split connected reductive group $G$.
The study of such representations is started by Barthel-Livn\'e~\cite{MR1290194,MR1361556} when $G = \GL_2(F)$.
They defined a notion of \emph{supersingular representations} and gave a classification of non-supersingular irreducible representations.
In particular, they proved that a representation is supersingular if and only if it is supercuspidal.
Here, a representation is called \emph{supercuspidal} if and only if it dose not appear as a subquotient of a parabolic induction from a proper parabolic subgroup.
By this theorem, to classify irreducible representations of $\GL_2(F)$, it is sufficient to classify irreducible supersingular representations.
When $G = \GL_2(\Q_p)$, irreducible supersingular representations are classified by Breuil~\cite{MR2018825}.
However, when $F\ne \Q_p$ a classification seems more complicated~\cite{breuil-Paskunas}.

Herzig~\cite{arXiv:1005.1713} gave a definition of a supersingular representation for any $G$ using the modulo $p$ Satake transform~\cite{arXiv:0910.4570}.
He also gave a classification of irreducible admissible representations in terms of supersingular representations when $G = \GL_n(F)$.
This is a generalization of a theorem of Barthel-Livn\'e.
In this paper, we generalize his classification to any $G$.

Now we state our main theorem.
Let $\overline{\kappa}$ be an algebraic closure of the residue field of $F$.
All representations in this paper are smooth representations over $\overline{\kappa}$.
Fix an $\mathcal{O}$-form of $G$ and denote it by the same letter $G$.
Let $K$ be the group of $\mathcal{O}$-valued points of $G$.
We also fix a Borel subgroup $B$ and a maximal torus $T\subset B$ of $G$.
Then we can define the notion of supersingular representation with respect to $(K,T,B)$. (See Herzig's paper~\cite[Definition~4.7]{arXiv:1005.1713} or Definition~\ref{defn:supersingular representation} in this paper.)
Let $\Pi$ be the set of simple roots.
Each subset $\Theta\subset\Pi$ corresponds to the standard parabolic subgroup $P_\Theta$.
Let $P_\Theta = M_\Theta N_\Theta$ be the Levi decomposition such that $T\subset M_\Theta$ and $N_\Theta$ is the unipotent radical of $P_\Theta$.
Consider the set $\mathcal{P}$ of all $\Lambda = (\Pi_1,\Pi_2,\sigma_1)$ such that:
\begin{itemize}
\item $\Pi_1$ and $\Pi_2$ are subsets of $\Pi$.
\item $\sigma_1$ is an irreducible admissible representation of $M_{\Pi_1}$ which is supersingular with respect to $(M_{\Pi_1}\cap K,T,M_{\Pi_1}\cap B)$.
\item Let $\omega_{\sigma_1}$ be the central character of $\sigma_1$ and put $\Pi_{\sigma_1} = \{\alpha\in\Pi\mid \langle\alpha,\check{\Pi}_1\rangle = 0,\ \omega_{\sigma_1}\circ\check{\alpha} = \trivrep_{\GL_1(F)}\}$.
Then $\Pi_2\subset\Pi_{\sigma_1}$.
\end{itemize}
Then the main theorem says that there exists a bijection between $\mathcal{P}$ and the set of isomorphism classes of irreducible admissible representations of $G$.

To state the theorem more precisely, we define the representation $I(\Lambda)$ for $\Lambda = (\Pi_1,\Pi_2,\sigma_1)\in \mathcal{P}$.
Let $P_\Lambda = M_\Lambda N_\Lambda$ be the Levi decomposition of the standard parabolic subgroup corresponding to $\Pi_1\cup \Pi_{\sigma_1}$.
First we construct the representation $\sigma_\Lambda$ of $M_\Lambda$.
Roughly, $\sigma_\Lambda$ is given by the following.
Since $\langle \Pi_1,\check{\Pi}_{\sigma_1}\rangle = 0$, $M_\Lambda$ is like the direct product of $M_{\Pi_1}$ and $M_{\Pi_{\sigma_1}}$.
We have a representation $\sigma_1$ of $M_{\Pi_1}$.
Since $\Pi_2\subset \Pi_{\sigma_1}$, $\Pi_2$ defines the standard parabolic subgroup of $M_{\Pi_{\sigma_1}}$ and it also defines the special representation~\cite{Grosse-special-rep}.
Let $\sigma_{\Lambda,2}$ be the special representation.
Then $\sigma_\Lambda$ is given by the tensor product $\sigma_1\boxtimes \sigma_{\Lambda,2}$.
However we have $M_\Lambda\not\simeq M_{\Pi_1}\times M_{\Pi_{\sigma_1}}$.
We define $\sigma_\Lambda$ as follows.

We can prove that $\sigma_1$ can be extended uniquely to $M_\Lambda$ such that $[M_{\Pi_{\sigma_1}},M_{\Pi_{\sigma_1}}]$ acts on it trivially (Lemma~\ref{lem:G/[M_2,M_2]simeq M_1/L_2}).
We denote the extended representation by the same letter $\sigma_1$.
Let $Q$ be the parabolic subgroup of $M_\Lambda$ corresponding to $\Pi_2\cup \Pi_{\sigma_1}$.
Then $Q$ defines the special representation of $M_\Lambda$~\cite{Grosse-special-rep} and denote it by $\sigma_{\Lambda,2}$.
From the definition of the special representation, the restriction of $\sigma_{\Lambda,2}$ to $M_{\Pi_{\sigma_1}}$ is the special representation of $M_{\Pi_{\sigma_1}}$ with respect to the standard parabolic subgroup corresponding to $\Pi_2$.
Now we define $\sigma_\Lambda = \sigma_1\otimes\sigma_{\Lambda,2}$ and put $I(\Lambda) = \Ind_{P_\Lambda}^G(\sigma_\Lambda)$.
The following is the main theorem of this paper.
\begin{thm}[Theorem~\ref{thm:Main theorem, general case}]
For $\Lambda\in \mathcal{P}$, $I(\Lambda)$ is irreducible and the correspondence $\Lambda\mapsto I(\Lambda)$ gives a bijection between $\mathcal{P}$ and the set of isomorphism classes of irreducible admissible representations.
\end{thm}
Using this theorem, we get the relation between supersingular representations and supercuspidal representations.
Recall that a representation is called \emph{supersingular} if it is supersingular with respect to any $(K,T,B)$.
\begin{thm}[Corollary~\ref{cor:supersinguality and supercuspidality}]
For an irreducible admissible representation $\pi$ of $G$, the following conditions are equivalent.
\begin{enumerate}
\item The representation $\pi$ is supersingular with respect to $(K,T,B)$.
\item The representation $\pi$ is supersingular.
\item The representation $\pi$ is supercuspidal.
\end{enumerate}
\end{thm}
We also give a criterion of the irreducibility of a principal series representation.
\begin{thm}
Let $\nu\colon T\to \overline{\kappa}^\times$ be a character.
Then $\Ind_B^G\nu$ is irreducible if and only if $\nu\circ\check{\alpha}\ne\trivrep_{\GL_1(F)}$ for all $\alpha\in\Pi$.
\end{thm}

Herzig proved much of his results for general $G$.
However, he proved the following two propositions only under some assumptions.
(Such assumptions are satisfied for $G = \GL_n$.)
\begin{enumerate}
\item A theorem of changing the weight~\cite[Corollary~6.10]{arXiv:1005.1713}.
\item A structure of a representation which have the ``trivial'' Satake parameter~\cite[Proposition~9.1]{arXiv:1005.1713}.
\end{enumerate}
In this paper, we prove these propositions for any $G$.
Then, Herzig's argument implies the main theorem.

We summarize the contents of this paper.
Using the modulo $p$ Satake transform, the notion of Satake parameters (or Hecke eigenvalues) is defined.
Such definition and properties are given in Section~\ref{sec:Satake parameters}.
A generalization of (1) and (2) is proved in Section~\ref{sec:A theorem of changing the weight}.
In Section~\ref{sec:A theorem of changing the weight}, we assume that the derived group of $G$ is simply connected.
Using these results, we prove our main theorem in Section~\ref{sec:Classification Theorem}.
Since we use results in Section~\ref{sec:A theorem of changing the weight}, first we assume that the derived group of $G$ is simply connected.
This assumption will be removed in subsection~\ref{subsec:General case and corollaries} using a $z$-extension.

\subsection*{Acknowledgments}
I thank Florian Herzig for reading the manuscript and giving helpful comments.
I also thank Tetsushi Ito for introducing me the theory of modulo $p$ representations of a $p$-adic group.

\section{Preliminaries}
\subsection{Notation}
In this paper, we use the following notation.
Let $p$ be a prime number, $F$ a finite extension of $\Q_p$, $\mathcal{O}$ its ring of integers, $\varpi \in \mathcal{O}$ a uniformizer, $\kappa = \mathcal{O}/(\varpi)$ the residue field and $q = \# \kappa$.
Let $G$ be a connected split reductive group over $\mathcal{O}$.
Fix a Borel subgroup $B\subset G$ and a split maximal torus $T\subset B$.
Let $U$ be the unipotent radical of $B$.
Then $B = TU$ is a Levi decomposition of $B$.
Let $\overline{B} = T\overline{U}$ be a Levi decomposition of the opposite group of $B$.
We also denote the group of $F$-valued points of $G$ by the same character $G$.
This should cause no confusion.
We use similar notation for other groups (for example, $B = B(F)$).
Set $K = G(\mathcal{O})$.
For any algebraic group $H$, let $H^\circ$ be the connected component containing the unit element and $Z_H$ the center of $H$.
For subgroups $H_1,H_2\subset H$, $Z_{H_1}(H_2) = \{h_1\in H_1\mid \text{$h_1h_2 = h_2h_1$ for all $h_2\in H_2$}\}$.
For a group $\Gamma$, $\trivrep_\Gamma$ is the trivial representation of $\Gamma$.
For a representation $V$ of $\Gamma$, $V^\Gamma$ is the space of invariants and $V_\Gamma$ is the space of coinvariants.

Let $(X^*,\Delta,X_*,\check{\Delta})$ be the root system of $(G,T)$.
Then $B$ determines the set of positive roots $\Delta^+\subset \Delta$ and the set of simple roots $\Pi\subset\Delta^+$.
Let $W$ be its Weyl group.
Let $\red\colon K = G(\mathcal{O})\to G(\kappa)$ be the canonical morphism.
The set of dominant (resp.~anti-dominant) elements in $X^*$ is denoted by $X^*_+$ (resp.~$X^*_-$).
We also use notation $X_{*,+}$ and $X_{*,-}$.
For $\lambda,\mu\in X_*$, we denote $\mu \le \lambda$ if $\lambda - \mu \in \Z_{\ge 0}\check{\Pi}$.

Let $P$ be a standard parabolic subgroup.
It has a Levi decomposition $P = MN$.
In this paper, we only consider the decomposition such that $T\subset M$.
The opposite parabolic subgroup of $P$ is denoted by $\overline{P} = M\overline{N}$.
We denote the Levi decomposition of the standard parabolic subgroup corresponding to $\Theta\subset\Pi$ by $P_\Theta = M_\Theta N_\Theta$.
The subset of $\Pi$ corresponding to $P$ is denoted by $\Pi_P$ or $\Pi_M$.
Put $\Delta_M = \Delta\cap \Z\Pi_M$ and $\Delta_M^+ = \Delta^+\cap \Delta_M$.
Let $W_M$ be the Weyl group of $\Delta_M$.
For $\nu\in X^*$, let $P_\nu = M_\nu N_\nu$ be the standard  parabolic subgroup corresponding to $\Pi_\nu = \{\alpha\in\Pi\mid \langle\nu,\check{\alpha}\rangle = 0\}$.
Put $W_\nu = \Stab_W(\nu)$, $\Delta_\nu = \{\alpha\in\Delta\mid \langle \nu,\check{\alpha}\rangle = 0\}$ and $\Delta_\nu^+ = \Delta^+\cap \Delta_\nu$.
(We will use these notation only when $\nu$ is dominant or anti-dominant. So the root system of $M_\nu$ is $\Delta_\nu$.)
We use similar notation for $\lambda\in X_*$.

For a subset $A\subset X^*$ and $A'\subset X_*$, $\langle A,A'\rangle = 0$ means $\langle \nu,\lambda\rangle = 0$ for all $\nu\in A$ and $\lambda\in A'$.
Notice that this condition is automatically satisfied if $A$ or $A'$ is empty.
We write $\langle A,\lambda\rangle = 0$ (resp.~$\langle \nu,A'\rangle = 0$) instead of $\langle A,\{\lambda\}\rangle = 0$ (resp.~$\langle \{\nu\},A'\rangle = 0$).

\subsection{Satake transform and irreducible representations of $K$}
Let $\overline{\kappa}$ be an algebraic closure of $\kappa$.
\emph{All representations in this paper are smooth representations over $\overline{\kappa}$}.
For a finite dimensional representation $V$ of $K$, let $\cInd_K^GV$ be a representation defined by
\[
	\cInd_K^GV = \{f\colon G\to V\mid f(xk) = k^{-1}f(x) \ (x\in G,k\in K),\ \text{$\supp f$ is compact}\}.
\]
The action of $g\in G$ is given by $(gf)(x) = f(g^{-1}x)$.
For $x\in G$ and $v\in V$, let $[x,v]\in \cInd_K^G(V)$ be an element defined by $\supp([x,v]) = xK$ and $[x,v](x) = v$.
Then $g[x,v] = [gx,v]$ and $[xk,v] = [x,kv]$ for $g\in G$ and $k\in K$.
For finite dimensional representations $V_1,V_2$ of $K$, $\Hom_G(\cInd_K^GV_1,\cInd_K^GV_2)$ is identified with
\begin{multline*}
	\mathcal{H}_G(V_1,V_2)\\
	= \left\{\varphi\colon G\to \Hom_{\overline{\kappa}}(V_1,V_2)\ \middle| \ 
	\begin{array}{l}
	\varphi(k_2xk_1) = k_2\varphi(x)k_1\ (k_1,k_2\in K,x\in G),\\
	\text{$\supp \varphi$ is compact}
	\end{array}
	\right\}.
\end{multline*}
The operator corresponding to $\varphi\in \mathcal{H}_G(V_1,V_2)$ is given by $f\mapsto \varphi * f$ where
\[
	(\varphi * f)(x) = \sum_{y\in G/K}\varphi(y)f(xy).
\]
We denote $\mathcal{H}_G(V,V)$ by $\mathcal{H}_G(V)$.
Let $\pi$ be a representation of $G$.
Then by the Frobenius reciprocity law, we have $\Hom_K(V,\pi)\simeq \Hom_G(\cInd_K^G(V),\pi)$.
Hence $\Hom_K(V,\pi)$ is a right $\mathcal{H}_G(V)$-module.
We denote the action of $\varphi\in\mathcal{H}_G(V)$ on $\psi\in \Hom_K(V,\pi)$ by $\psi*\varphi$.

When $V$ is irreducible, the structure of $\mathcal{H}_G(V)$ is given by the \emph{Satake transform}~\cite{arXiv:0910.4570}.
Namely, the Satake transform $S_G\colon \mathcal{H}_G(V)\to \mathcal{H}_T(V^{\overline{U}(\kappa)})$ defined by
\[
	S_G(\varphi)(t) = \sum_{u\in \overline{U}/\overline{U}(\mathcal{O})}\varphi(tu)|_{V^{\overline{U}(\kappa)}}
\]
is injective and its image is $\{\varphi\in \mathcal{H}_T(V^{\overline{U}(\kappa)})\mid \supp\varphi\subset T_+\}$ where $T_+ = \{t\in T\mid \alpha(t)\in \mathcal{O}\ (\alpha\in\Delta^+)\}$.
A homomorphism $X_*\times T(\mathcal{O})\to T$ defined by $(\lambda,t_0)\mapsto \lambda(\varpi)t_0$ is an isomorphism and it induces $X_{*,+}\times T(\mathcal{O})\simeq T_+$.
Hence $S_G$ gives an isomorphism $\mathcal{H}_G(V)\simeq\overline{\kappa}[X_{*,+}]$.
For $\lambda\in X_{*,+}$, there exists $T_\lambda\in \mathcal{H}_G(V)$ such that $\supp T_\lambda = K\lambda(\varpi)K$ and $T_\lambda(\lambda(\varpi))$ is given by $V\twoheadrightarrow V_{\overline{N_\lambda}(\kappa)}\simeq V^{N_\lambda(\kappa)}\hookrightarrow V$.
Then $\{T_\lambda\mid \lambda\in X_{*,+}\}$ gives a basis of $\mathcal{H}_G(V)$.
When we want to emphasis the group $G$, we write $T^G_\lambda$ instead of $T_\lambda$.
For $\lambda\in X_{*}$, let $\tau_\lambda\in \overline{\kappa}[X_{*}]$ be an element corresponding to $\lambda$.
Then $\{\tau_\lambda\mid \lambda\in X_{*,+}\}$ gives a basis of $\overline{\kappa}[X_{*,+}]$.
The relation between $S_G(T_\lambda)$ and $\tau_\lambda$ is given by Herzig~\cite[Proposition~5.1]{arXiv:1005.1713}.
An algebra homomorphism $\overline{\kappa}[X_{*,+}]\to \overline{\kappa}$ is parameterized by $(M,\chi_M)$ where $M$ is the 
Levi subgroup of a standard parabolic subgroup and $\chi_M$ is a group homomorphism $X_{M,*,0}\to \overline{\kappa}^\times$ where $X_{M,*,0} = \{\lambda\in X_*\mid \langle \lambda,\Pi_M\rangle = 0\}$~\cite[Proposition~4.1]{arXiv:1005.1713}.
Therefore, an algebra homomorphism $\mathcal{H}_G(V)\to \overline{\kappa}$ is parameterized by the same pair.
\begin{rem}\label{rem:choice of uniformizer and Hecke algebra}
Since an isomorphism $\mathcal{H}_T(V^{\overline{U}(\kappa)})\simeq\overline{\kappa}[X_*]$ depends on a choice of a uniformizer $\varpi$, the above parameterization is not natural.
More natural way is given by Herzig~\cite{arXiv:0910.4570}.
In this paper, we fix a uniformizer and identify $\mathcal{H}_G(V)$ with $\overline{\kappa}[X_{*,+}]$.
(It is only for a simplification of notation.)
\end{rem}

Let $P = MN$ be the Levi decomposition of a standard parabolic subgroup.
Then the partial Satake transform $S_G^M\colon \mathcal{H}_G(V)\to \mathcal{H}_M(V^{\overline{N}(\kappa)})$ is injective and it satisfies $S_M\circ S_G^M = S_G$~\cite[2.3]{arXiv:1005.1713}.
Assume that $\chi\colon\mathcal{H}_G(V)\to \overline{\kappa}$ is parameterized by $(M,\chi_M)$.
Then $M$ is characterized by the following property: $\chi$ factors through $S_G^{M'}$ if and only if $M'\supset M$.
We also have the following: $\chi_M(\lambda) = \chi(\tau_\lambda)$.

Let $V_1,V_2$ be irreducible representations of $K$.
For each $\lambda\in X_{*,+}$, there exists $\varphi\in \mathcal{H}(V_1,V_2)\setminus\{0\}$ whose support is $K\lambda(\varpi)K$ if and only if $V_1^{\overline{N_\lambda}(\kappa)}\simeq V_2^{\overline{N_\lambda}(\kappa)}$ as $M_\lambda(\kappa)$-representations.
Moreover such $\varphi$ is unique up to constant multiple.
The homomorphism $\varphi(\lambda(\varpi))$ is given by $V_1\twoheadrightarrow (V_1)_{N_\lambda(\kappa)}\simeq V_2^{\overline{N_\lambda}(\kappa)}\hookrightarrow V_2$.
(See the proof of \cite[Proposition~6.3]{arXiv:1005.1713}.)

All irreducible representations of $K$ factor through $K\to G(\kappa)$.
If the derived group of $G$ is simply connected, such representation is parameterized by its lowest weight.
If $\nu\in X^*$ satisfies $-q <  \langle \nu,\check{\alpha}\rangle \le 0$ for all $\alpha\in\Pi$ then the restriction of the irreducible representation of $G(\overline{\kappa})$ with lowest weight $\nu$ to $G(\kappa)$ is irreducible and they give all irreducible representations of $G(\kappa)$.
When $V$ is the restriction of an irreducible representation with lowest weight $\nu$, we call $\nu$ a lowest weight of $V$.
(For $\nu_0\in X^*$ such that $\langle \nu_0,\check{\Pi}\rangle = 0$, the restriction of the irreducible representations with lowest weight $\nu$ and $\nu + (q - 1)\nu_0$ are isomorphic to each other.
Hence $\nu$ is not determined by $V$ uniquely.)

\section{Satake parameters}\label{sec:Satake parameters}
\subsection{Definition and some lemmas}\label{subsec:stake parameters:definition and some lemmas}
We start with the following definition.
\begin{defn}\label{defn:Satake parameters}
Let $\pi$ be a representation of $G$.
An algebra homomorphism $\chi\colon \overline{\kappa}[X_{*,+}]\to\overline{\kappa}$ is called a \emph{Satake parameter of $\pi$} if there exist an irreducible $K$-representation $V$ and $\psi\in\Hom_K(V,\pi)\setminus\{0\}$ such that for all $\varphi\in\mathcal{H}_G(V)$, $\psi*\varphi = \chi(S_G(\varphi))\psi$.
\end{defn}
Let $\Satakepar(\pi,V)$ be the set of Satake parameters appearing in $\Hom_K(V,\pi)$.
We denote the set of Satake parameters of $\pi$ by $\Satakepar(\pi)$.
Then we have $\Satakepar(\pi) = \bigcup_V\Satakepar(\pi,V)$.
If $\pi$ is admissible, then $\Satakepar(\pi)\ne \emptyset$.
We give some propositions about Satake parameters.
Before proving some properties of Satake parameters, we give some fundamental facts about a structure of $G$.
\begin{lem}\label{lem:G/[M_2,M_2]simeq M_1/L_2}
Let $\Pi = \Pi_1\cup \Pi_2$ be a partition of $\Pi$ such that $\langle \Pi_1,\check{\Pi}_2\rangle = 0$ and $P_i =M_iN_i$ the standard parabolic subgroup corresponding to $\Pi_i$.
Let $L_2$ be the subgroup of $T\subset M_1$ generated by $\{\Imm\check{\alpha}\mid \alpha\in\Pi_2\}$.
Then we have $G/[M_2,M_2]\simeq M_1/L_2$.
\end{lem}
\begin{proof}
Let $\overline{F}$ be a separable closure of $F$.
In this proof, we write $\mathbf{G} = G(\overline{F})$.
(The same notation is used for other groups.)
Set $\check{\Pi}_2^\perp = \{\nu\in X^*\mid \langle \nu,\check{\Pi}_2\rangle = 0\}$.
Since $\mathbf{G}/[\mathbf{M}_2,\mathbf{M}_2]$ and $\mathbf{M}_1/\mathbf{L}_2$ have the same root data $(\check{\Pi}_2^\perp,\Delta_{M_1},X_*/(\Q\check{\Pi}_2\cap X_*),\check{\Delta}_{M_1})$, these are isomorphic.
For an abelian group $A$, let $A_{\mathrm{tors}}$ be its torsion part.
By a theorem of Kottwitz, the Galois cohomology $H^1(F,[\mathbf{M}_2,\mathbf{M}_2])$ (resp.~$H^1(F,\mathbf{G})$) is isomorphic to $((X_*\cap \Q\check{\Pi}_2)/\Z\check{\Pi}_2)_{\mathrm{tors}}$ (resp.~$(X_*/\Z\check{\Pi})_{\mathrm{tors}}$).
Hence $H^1(F,[\mathbf{M}_2,\mathbf{M}_2])\to H^1(F,\mathbf{G})$ is injective.
Therefore, $(\mathbf{G}/[\mathbf{M}_2,\mathbf{M}_2])^{\Gal(\overline{F}/F)} = G/[M_2,M_2]$.
Since $\mathbf{L}_2$ is a torus, $H^1(F,\mathbf{L}_2)$ is trivial.
Hence $(\mathbf{M}_1/\mathbf{L}_2)^{\Gal(\overline{F}/F)} = M_1/L_2$.
The lemma follows.
\end{proof}
\begin{prop}
There is a one-to-one correspondence between the character $\nu_G$ of $G$ and the character $\nu_T$ of $T$ such that $\nu_T\circ\check{\alpha}$ is trivial for all $\alpha\in\Pi$.
It is characterized by $\nu_T = \nu_G|_T$.
\end{prop}
\begin{proof}
Apply the previous lemma for $\Pi_1 = \emptyset$ and $\Pi_2 = \Pi$.
\end{proof}
\begin{cor}\label{cor:char of K and G}
Let $\nu_K$ be a character of $K$.
Then there exists a character $\nu_G$ of $G$ such that $\nu_K = \nu_G|_K$.
If $\nu_G(\lambda(\varpi)) = 1$ for all $\lambda\in X_*$, then $\nu_G$ is unique.
\end{cor}
\begin{proof}
If the derived group of $G$ is simply connected, it is known that $\nu_K$ has a lowest weight $\nu$ which satisfies $(\nu\circ \check{\alpha})(\mathcal{O}^\times) = 1$ for all $\alpha\in \Pi$.
Therefore, the corollary follows from the above proposition.
In general, let $1\to Z\to \widetilde{G}\to G\to 1$ be a $z$-extension of $G$, $\widetilde{K}$ the group of $\mathcal{O}$-valued points of $\widetilde{G}$ and $\widetilde{T}$ the inverse image of $T$ in $\widetilde{G}$.
Then there exists a character $\nu_{\widetilde{G}}$ such that $\nu_{\widetilde{G}}|_{\widetilde{K}}$ is a pull-back of $\nu_K$ and $\nu_{\widetilde{G}}(\lambda(\varpi)) = 1$ for all $\lambda\in X_*(\widetilde{T})$.
Hence $\nu_{\widetilde{G}}|_Z$ is trivial.
Therefore, it gives a character $\nu_G$ of $G$ and $\nu_G|_K = \nu_K$.
\end{proof}
For a character $\nu$ of $G$, $\varphi\mapsto (g\mapsto \varphi_\nu(g) = \varphi(g)\nu(g))$ gives an isomorphism $\mathcal{H}_G(V)\simeq \mathcal{H}_G(V\otimes\nu|_K)$.
\begin{lem}\label{lem:Satake and character}
For a parabolic subgroup $P = MN$, the homomorphism $\varphi\mapsto \varphi_\nu$ is compatible with the partial Satake transform $S_G^M$.
\end{lem}
\begin{proof}
We have
\[
	(S_G^M\varphi_\nu)(m) = \sum_{n\in N/(N\cap K)}\nu(mn)\varphi(mn).
\]
Since $N\subset [G,G]$, we have $\nu(n) = 1$.
Therefore,
\[
	\sum_{n\in N/(N\cap K)}\nu(mn)\varphi(mn)
	=
	\nu(m)\sum_{n\in N/(N\cap K)}\varphi(mn)
	=
	\nu(m)(S_G^M\varphi)(m).
\]
\end{proof}

Now we give some properties on Satake parameters.
The following proposition is obvious.
\begin{prop}\label{prop:Satake parameters of subrepresentation}
If $\pi'\subset\pi$, then $\Satakepar(\pi',V)\subset\Satakepar(\pi,V)$.
\end{prop}

The following proposition follows from \cite[Lemma~2.14]{arXiv:1005.1713}.
\begin{prop}\label{prop:Satake parameters of parbolic induction}
Let $P = MN$ be a parabolic subgroup, $\sigma$ a representation of $M$ and $V$ an irreducible representation of $K$.
Then we have $\Satakepar(\Ind_P^G(\sigma),V) = \Satakepar(\sigma,V^{\overline{N}(\kappa)})|_{\overline{\kappa}[X_{*,+}]}$.
In particular, we have $\Satakepar(\Ind_P^G(\sigma)) = \Satakepar(\sigma)|_{\overline{\kappa}[X_{*,+}]}$.
\end{prop}

Let $\chi_1,\chi_2\colon \overline{\kappa}[X_{*,+}]\to\overline{\kappa}$ be algebra homomorphisms.
Define $\chi_1\otimes\chi_2\colon \overline{\kappa}[X_{*,+}]\to \overline{\kappa}$ by $(\chi_1\otimes\chi_2)(\tau_\lambda) = \chi_1(\tau_\lambda)\chi_2(\tau_\lambda)$.
\begin{prop}
Assume $\chi_i$ is parameterized by $(M_i,\chi_{M_i})$.
Then $\chi_1\otimes\chi_2$ is parameterized by $(M,\chi_M)$ where $\Pi_M = \Pi_{M_1}\cap \Pi_{M_2}$ and $\chi_M = \chi_{M_1}|_{X_{M,*,0}}\chi_{M_2}|_{X_{M,*,0}}$.
\end{prop}
\begin{proof}
This follows from the argument of the proof of \cite[Corollary~1.5]{arXiv:0910.4570}.
\end{proof}

\begin{prop}
Let $\nu$ be a character of $G$ and $\pi$ a representation of $G$.
Then $\Satakepar(\pi\otimes\nu) = \Satakepar(\pi)\otimes \chi_\nu$ here $\chi_\nu\colon \overline{\kappa}[X_{*,+}]\to \overline{\kappa}$ is given by $\chi_\nu(\tau_\lambda) = \nu(\lambda(\varpi))$.
\end{prop}
\begin{proof}
This follows from Lemma~\ref{lem:Satake and character}.
\end{proof}

\begin{prop}\label{prop:Satake parameter of a character}
Let $\nu$ be a character of $G$.
Then $\Satakepar(\nu) = \{\chi_\nu\}$.
\end{prop}
\begin{proof}
We have an injective homomorphism $\nu\hookrightarrow \Ind_B^G(\nu|_T)$.
Hence we have $\emptyset\ne\Satakepar(\nu)\subset \Satakepar(\Ind_B^G(\nu|_T)) = \Satakepar(\nu|_T)|_{\overline{\kappa}[X_{*,+}]} = \{\chi_\nu\}$.
\end{proof}

\subsection{Restriction and Satake parameter}
Let $G_1$ be a connected subgroup of $G$ which contains $[G,G]$.
Let $X_{G_1,*}$ be the group of cocharacters of $G_1\cap T$.
Put $X_{G_1,*,+} = X_{*,+}\cap X_{G_1,*}$.
Then we have $\mathcal{H}_{G_1}(V)\simeq\overline{\kappa}[X_{G_1,*,+}]$.
Since $X_{G_1,*,+}\subset X_{*,+}$, we have an injective homomorphism $\overline{\kappa}[X_{G_1,*,+}]\hookrightarrow \overline{\kappa}[X_{*,+}]$.
This induces $\Phi\colon \mathcal{H}_{G_1}(V)\hookrightarrow \mathcal{H}_G(V)$.
\begin{lem}\label{lem:hecke algebra of G_1}
We have $\Imm\Phi = \{\varphi\in \mathcal{H}_G(V)\mid \supp\varphi\subset G_1K\}$ and an isomorphism $\Imm\Phi\simeq \mathcal{H}_{G_1}(V)$ is given by $\varphi\mapsto \varphi|_{G_1}$.
\end{lem}
\begin{proof}
Put $\mathcal{H}_1 = \{\varphi\in \mathcal{H}_G(V)\mid \supp \varphi\subset G_1K\}$.
Then $\mathcal{H}_1$ has a basis $\{T^G_\lambda\mid \lambda\in X_{G_1,*,+}\}$.
To prove the first statement of the lemma, it is sufficient to prove that if $\lambda\in X_{G_1,*,+}$ then $S_G(T^G_\lambda)\in \overline{\kappa}[X_{G_1,*,+}]$ and $\{S_G(T^G_\lambda)\mid \lambda\in X_{G_1,*,+}\}$ is a basis of $\overline{\kappa}[X_{G_1,*,+}]$.
We have $S_G(T^G_\lambda)\in \tau_\lambda + \sum_{\mu\le \lambda}\overline{\kappa}\tau_\mu$.
Since $\check{\Pi}\subset X_{G_1,*}$, $\lambda\in X_{G_1,*}$ and $\mu \le \lambda$ imply $\mu \in X_{G_1,*}$.
Therefore we get the first statement.

Since $U$ is the unipotent radical of the Borel subgroup $B\cap G_1$ of $G_1$, we have $S_G(T_\lambda^G) = S_{G_1}(T_\lambda^G|_{G_1})$ for $\lambda\in X_{G_1,*,+}$ by the definition of the Satake transform.
We get the second statement.
\end{proof}
\begin{prop}\label{prop:comparison of Satake parameters, derived group}
Let $\pi$ be a representation of $G$ and $V$ an irreducible representation of $K$.
Then we have $\Satakepar(\pi,V)|_{\overline{\kappa}[X_{G_1,*,+}]}\subset \Satakepar(\pi|_{G_1},V|_{G_1\cap K})$.
Hence $\Satakepar(\pi)|_{\overline{\kappa}[X_{G_1,*,+}]}\subset\Satakepar(\pi|_{G_1})$.

Moreover, if $\pi$ has a central character, then for each irreducible $(G_1\cap K)$-representation $V_1$, we have $\Satakepar(\pi|_{G_1},V_1) = \bigcup_{V|_{G_1\cap K} = V_1}\Satakepar(\pi,V)|_{\overline{\kappa}[X_{G_1,*,+}]}$.
Hence $\Satakepar(\pi|_{G_1}) = \Satakepar(\pi)|_{\overline{\kappa}[X_{G_1,*,+}]}$.
\end{prop}
\begin{proof}
Put $K_1 = G_1\cap K$.
Since $[K,K]\subset K_1$, the restriction of an irreducible $K$-representation to $K_1$ is also irreducible.
Let $V$ be an irreducible representation of $K$.
We prove $\Satakepar(\pi,V)|_{\overline{\kappa}[X_{G_1,*,+}]}\subset\Satakepar(\pi|_{G_1},V|_{K_1})$.
It is sufficient to prove that
\[
	\Hom_K(V,\pi)\hookrightarrow \Hom_{K_1}(V,\pi)
\]
is $\mathcal{H}_{G_1}(V)$-module homomorphism.
Let $\varphi\in\mathcal{H}_{G_1}(V)$ and $\psi\in\Hom_K(V,\pi)$.
Then for each $v\in V$, 
\[
	(\psi * \Phi(\varphi))(v) = \sum_{g\in G/K}g\psi(\Phi(\varphi)(g^{-1})v)
	= \sum_{g\in G_1K/K}g\psi(\Phi(\varphi)(g^{-1})v).
\]
The claim follows from $G_1/K_1\simeq G_1K/K$.

Assume that $\pi$ has a central character.
Let $V_1$ be an irreducible representation of $K_1$.
Fix an irreducible representation $V$ of $K$ such that $V|_{K_1} = V_1$.
Take such representation such that a central character of $V$ is the same as that of $\pi$.
Set $K' = K_1Z_K$.
Then $K'$ has a finite index in $K$ and we have
\[
	\Hom_{K_1}(V,\pi) = \Hom_{K'}(V,\pi)\simeq\Hom_K(\Ind_{K'}^K(V),\pi).
\]
Since $V$ has a structure of a representation of $K$, we have $\Ind_{K'}^K(V)\simeq\Ind_{K'}^K(\trivrep_{K'})\otimes V$.
Therefore we have
\[
	\Psi\colon \Hom_{K_1}(V,\pi)\simeq\Hom_K(\Ind_{K'}^K(\trivrep_{K'})\otimes V,\pi).
\]
Explicitly, this isomorphism is given by
\[
	\Psi(\psi)(f\otimes v) = \sum_{x\in K/K'}f(x)x\psi(x^{-1}(v)).
\]
Therefore, for $\varphi\in\mathcal{H}_{G_1}(V)$, we have
\begin{multline*}
	\Psi(\psi *\varphi)(f\otimes v) = \sum_{x\in K/K'}f(x)x\sum_{g\in G_1/K_1}g\psi(\varphi(g^{-1})x^{-1}v)\\
	= \sum_{x\in K/K'}\sum_{g\in G_1/K_1}f(x)(xg)\psi(\varphi((xg)^{-1})v).
\end{multline*}
Replacing $g$ with $x^{-1}gx$, we have
\begin{multline*}
	\Psi(\psi *\varphi)(f\otimes v) 
	= \sum_{x\in K/K'}\sum_{g\in G_1/K_1}f(x)gx\psi(x^{-1}\varphi(g^{-1})v)\\
	= \sum_{g\in G_1/K_1}g\Psi(\psi)(f\otimes \varphi(g^{-1})v).
\end{multline*}
Since $K'$ is a normal subgroup of $K$ and $K/K'$ is commutative, the representation $\Ind_{K'}^K(\trivrep_{K'})$ has a filtration $\{X_i\}$ such that $X_i/X_{i - 1}\simeq \nu_i$ for some character $\nu_i$ of $K$.
Set $X = \Ind_{K'}^K(\trivrep_{K'})\otimes V$, $Y = \Hom_K(X,\pi)$ and $Y_i = \Hom_K(X/X_i,\pi)$.
Then $\{Y_i\}$ is a filtration of $Y$ and $Y_i/Y_{i - 1}\hookrightarrow \Hom_K(\nu_i\otimes V,\pi)$.
By the above formula, $Y_i$ is stable under the action of $\varphi\in \mathcal{H}_{G_1}(V)$.
Hence $\varphi$ acts on $Y_i/Y_{i - 1}$.
Extend $\nu_i$ to a character of $G$ such that $\nu_i$ is trivial on $G_1$.
Then we have $\mathcal{H}_G(V)\simeq \mathcal{H}_G(\nu_i\otimes V)$ by $\varphi'\mapsto \varphi'_{\nu_i}$.
We have an action of $\Phi(\varphi)_{\nu_i}\in \mathcal{H}_G(\nu_i\otimes V)$ on $\Hom_K(\nu_i\otimes V,\pi)$.
We prove that these actions are compatible with $Y_i/Y_{i - 1}\hookrightarrow \Hom_K(\nu_i\otimes V,\pi)$.

Since $\nu_i$ is trivial on $G_1$, we have $a\otimes\varphi(g^{-1})v = \Phi(\varphi)_{\nu_i}(g^{-1})(a\otimes v)$ for $g\in G_1$.
The function $g\mapsto g\Psi(\psi)(\Phi(\varphi)_{\nu_i}(g^{-1})(a\otimes v))$ is right $K$-invariant.
Therefore, 
\begin{multline*}
	\sum_{g\in G_1/K_1}g\Psi(\psi)(a\otimes\varphi(g^{-1})v)
	=
	\sum_{g\in G_1K/K}g\Psi(\psi)(\Phi(\varphi)_{\nu_i}(g^{-1})(a\otimes v))\\
	=
	\sum_{g\in G/K}g\Psi(\psi)(\Phi(\varphi)_{\nu_i}(g^{-1})(a\otimes v))
	=
	(\Psi(\psi) *\Phi(\varphi)_{\nu_i})(a\otimes v).
\end{multline*}
This means that the actions are compatible.

Hence each element of $\Satakepar(\pi|_{G_1},V)$ appears in $\Satakepar(\pi,\nu_i\otimes V)|_{\overline{\kappa}[X_{G_1,*,+}]}$ for some $i$.
Since $\nu_i$ is trivial on $K_1$, $(\nu_i\otimes V)|_{K_1}\simeq V|_{K_1}\simeq V_1$.
We get $\Satakepar(\pi|_{G_1},V)\subset\bigcup_{V'|_{K_1} = V|_{K_1}}\Satakepar(\pi,V')|_{\overline{\kappa}[X_{G_1,*,+}]}$.
\end{proof}

\subsection{Satake parameter of tensor product}
Consider the setting in Lemma~\ref{lem:G/[M_2,M_2]simeq M_1/L_2}.
Namely, let $\Pi = \Pi_1\cup\Pi_2$ be a partition of $\Pi$ such that $\langle \Pi_1,\check{\Pi}_2\rangle = 0$.
Let $P_i = M_iN_i$ be the standard parabolic subgroup corresponding to $\Pi_i$.
Set $H_1 = Z_{M_1}([M_2,M_2])^\circ$.
Put $\Pi_2^\perp = \{\lambda\in X_*\mid \langle \lambda,\Pi_2\rangle = 0\}$.
Then the group of cocharacters of $H_1\cap T$ is $\Pi_2^\perp$.
We also have $[M_1,M_1]\subset H_1\subset M_1$.
Put $X_{H_1,*,+} = X_{*,+}\cap \Pi_2^\perp$.

\begin{lem}
If $\lambda,\mu\in X_{*,+}$ satisfies $\mu\le \lambda$ and $\lambda\in X_{H_1,*,+}$, then $\lambda - \mu\in\Z_{\ge 0}\Pi_1$.
In particular, $\mu\in X_{H_1,*,+}$.
\end{lem}
\begin{proof}
For each $\alpha\in\Pi$, take $n_\alpha\in\Z_{\ge 0}$ such that $\lambda - \mu = \sum_{\alpha\in\Pi}n_\alpha \check{\alpha}$.
Then for $\beta\in\Pi_2$, we have $\sum_{\alpha\in\Pi_2}n_\alpha\langle\beta,\check{\alpha}\rangle = -\langle \beta,\mu\rangle \le 0$.
Since $(\langle\beta,\check{\alpha}\rangle)_{\alpha,\beta\in\Pi_2}$ is positive definite, we have $n_\alpha = 0$ for all $\alpha\in\Pi_2$.
\end{proof}

Fix an irreducible representation $V$ of $K$ and put $V_1 = V^{\overline{N}_1(\kappa)}$.
Then $V_1$ is irreducible as a representation of $M_1\cap K$.
Since $[M_1,M_1]\subset H_1\subset M_1$, $V_1$ is also irreducible as a representation of $H_1\cap K$.
We have $\overline{\kappa}[X_{H_1,*,+}]\hookrightarrow \overline{\kappa}[X_{*,+}]$.
Hence we get $\Phi'\colon \mathcal{H}_{H_1}(V_1)\hookrightarrow \mathcal{H}_G(V)$.
By the above lemma and the argument in the proof of Lemma~\ref{lem:hecke algebra of G_1}, we get the following lemma.

\begin{lem}
We have $\Imm\Phi' = \{\varphi\in \mathcal{H}_G(V)\mid \supp\varphi\subset KH_1K\}$ and the isomorphism $\Imm\Phi'\simeq \mathcal{H}_{H_1}(V)$ is given by $\varphi\mapsto \varphi|_{H_1}$.
\end{lem}

Let $\pi$ be a representation of $G$.
Consider the following homomorphism
\[
	\Hom_K(V,\pi)\to \Hom_{M_1\cap K}(V_1,\pi).
\]
Since $V$ is generated by $V_1$ as a $K$-representation, this is injective.
The left hand side is $\mathcal{H}_G(V)\simeq\overline{\kappa}[X_{*,+}]$-module and the right hand side is $\mathcal{H}_{M_1}(V_1)\simeq\overline{\kappa}[X_{M_1,*,+}]$-module where $X_{M_1,*,+} = \{\lambda\in X_*\mid \langle \lambda,\alpha\rangle\ge 0\  (\alpha\in\Pi_{M_1})\}$.
Therefore, both sides are $\overline{\kappa}[X_{H_1,*,+}]$-modules.
We prove that the above embedding is a $\overline{\kappa}[X_{H_1,*,+}]$-modules homomorphism.
We need a lemma.
\begin{lem}
For $m\in M_1$ and $n\in N_1$, if $mn\in KH_1K$, then $n\in K$.
\end{lem}
\begin{proof}
By a Cartan decomposition, we can choose $\lambda\in X_{H_1,*,+}$, $\lambda_1\in X_{M_1,*,+}$ and $k_1\in M_1\cap K$ such that $mn\in K\lambda(\varpi)K$ and $m\in (M_1\cap K)\lambda_1(\varpi)k_1$.
Then we have $\lambda_1(\varpi)(k_1nk_1^{-1})\in K\lambda(\varpi)K$.
Put $n_1 = k_1nk_1^{-1}\in N_1$.
We prove $n_1\in K$.

By the assumption, we have $N_1\subset M_2$.
Therefore, $\lambda_1(\varpi)n_1$ is in $M_2$.
Take $\lambda_2\in X_{M_2,*,+}$ such that $\lambda_1(\varpi)n_1\in (M_2\cap K)\lambda_2(\varpi)(M_2\cap K)$.
Then $K\lambda_2(\varpi)K\cap K\lambda(\varpi)K\ne\emptyset$.
Therefore, $\lambda_2\in W\lambda$.
The Weyl group $W$ preserves each connected component of a root system.
Hence $W$ preserves $\check{\Pi}_2^\perp$.
Hence $\lambda_2\in \Pi_2^\perp$.
Therefore, $\lambda_2(\varpi)$ commutes with an element of $M_2$.
Hence $\lambda_1(\varpi)n_1\in (M_2\cap K)\lambda_2(\varpi)(M_2\cap K) = \lambda_2(\varpi)(M_2\cap K)$.
Therefore, $\lambda_2(\varpi)^{-1}\lambda_1(\varpi)n_1\in K$.
We get $n_1\in K$.
\end{proof}

\begin{lem}\label{lem:comparison of satake parameters, restriction to Levi}
Let $\pi$ be a representation of $G$.
The homomorphism
\[
	\Hom_K(V,\pi)\to \Hom_{M_1\cap K}(V_1,\pi).
\]
is a $\overline{\kappa}[X_{H_1,*,+}]$-module homomorphism.
\end{lem}
\begin{proof}
Let $\varphi\in \mathcal{H}_{H_1}(V_1)$.
Take $\psi\in \Hom_K(V,\pi)$ and $v\in V_1$.
We have
\begin{multline*}
(\psi * \Phi'(\varphi))(v) = \sum_{g\in G/K}g\psi(\Phi'(\varphi)(g^{-1})v)\\
= \sum_{m\in M_1/(M_1\cap K)}\sum_{n\in N_1/(N_1\cap K)}mn\psi(\Phi'(\varphi)(n^{-1}m^{-1})v)
\end{multline*}
Since $\supp\Phi'(\varphi)\subset KH_1K$, $\Phi'(\varphi)(n^{-1}m^{-1}) = 0$ if $n\not\in N_1\cap K$ by the above lemma.
Therefore, we have
\begin{multline*}
(\psi *\Phi'(\varphi))(v) = \sum_{m\in M_1/(M_1\cap K)}m\psi(\Phi'(\varphi)(m^{-1})v)\\ = \sum_{m\in M_1/(M_1\cap K)}m\psi(\Phi'(\varphi)(m^{-1})v).
\end{multline*}
By the same argument, we have $S_G^{M_1}(\Phi'(\varphi))(m) = \Phi'(\varphi)(m)$ for $m\in M_1$.
Therefore, we get the lemma.
\end{proof}

Let $\pi_1,\pi_2$ be representations of $G$ with the central characters such that $[M_1,M_1]$ acts on $\pi_2$ trivially and the center of $M_2$ acts on $\pi_2$ as a scalar.
Put $\pi = \pi_1\otimes\pi_2$.

\begin{rem}
The group $H_1$ is generated by $H_1\cap T$ and one-dimensional unipotent subgroup corresponding to each $\alpha\in\Delta\cap\Z\Pi_1$.
Since $H_1\cap T\subset Z_{M_2}^\circ$ and one-dimensional unipotent subgroup corresponding to $\alpha\in\Delta\cap\Z\Pi_1$ is a subgroup of $[M_1,M_1]$, $H_1$ is generated by $[M_1,M_1]$ and $Z_{M_2}^\circ$.
Therefore, $H_1$ acts on $\pi_2$ by a scalar.
\end{rem}

\begin{prop}\label{prop:Satake parameters of tensor product}
We have $\Satakepar(\pi)|_{\overline{\kappa}[X_{H_1,*,+}]} \subset \Satakepar(\pi_1|_{H_1})\otimes \Satakepar(\pi_2|_{H_1})$.
\end{prop}
\begin{proof}
We have $\Satakepar(\pi)|_{\overline{\kappa}[X_{H_1,*,+}]}\subset\Satakepar(\pi|_{M_1})|_{\overline{\kappa}[X_{H_1,*,+}]}$ by the above lemma.
By Proposition~\ref{prop:comparison of Satake parameters, derived group}, we have $\Satakepar(\pi|_{M_1})|_{\overline{\kappa}[X_{H_1,*,+}]} \subset \Satakepar(\pi|_{H_1})$.
Since $H_1$ acts on $\pi_2$ by a scalar, $\Satakepar(\pi|_{H_1})= \Satakepar(\pi_1|_{H_1})\otimes\Satakepar(\pi_2|_{H_1})$ by Lemma~\ref{lem:Satake and character}.
\end{proof}

\begin{lem}
If $\pi_1$ is admissible as a representation of $[M_1,M_1]$ and $\pi_2$ is admissible as a representation of $G$.
Then $\pi$ is admissible as a representation of $G$.
\end{lem}
\begin{proof}
Let $K'$ be a compact open subgroup.
Then we have $\pi^{K'} = (\pi_1^{[M_1,M_1]\cap K'}\otimes \pi_2)^{K'}$.
Since $\pi_1^{[M_1,M_1]\cap K'}$ is finite dimensional and $\pi_2$ is admissible, $\dim\pi^{K'} < \infty$.
\end{proof}

We give some corollaries of Proposition~\ref{prop:Satake parameters of tensor product} which we will use.
We suppose the following additional assumptions.
\begin{itemize}
\item The representation $\pi_1$ is an admissible $[M_1,M_1]$-representation.
\item The representation $\pi_2$ is an admissible $G$-representation.
\item The derived group $[M_2,M_2]$ acts on $\pi_1$ trivially and the center of $M_1$ acts on $\pi_1$ as a scalar.
\item We have $\#\Satakepar(\pi_1|_{M_1}) = \#\Satakepar(\pi_2|_{M_2}) = 1$.
\end{itemize}
Then there exists a unique parabolic subgroup $P = MN$ such that $\Satakepar(\pi_1|_{M_1}) = \{\chi_1 = (M\cap M_1,\chi_{M_1})\}$ and $\Satakepar(\pi_2|_{M_2}) = \{\chi_2 = (M\cap M_2,\chi_{M_2})\}$ for some $\chi_{M_1}$ and $\chi_{M_2}$.
By the above lemma, $\pi$ is admissible.

\begin{cor}\label{cor:zero-point of Satake parameter of tensor}
Any $\chi\in\Satakepar(\pi)$ is parameterized by $(M,\chi_M)$ for some $\chi_M$.
\end{cor}
\begin{proof}
Take $M'$ such that $\chi$ is parameterized by $(M',\chi_{M'})$.
For each $\alpha\in\Pi$, take $\lambda_\alpha\in X_{*,+}$ such that $\langle \Pi\setminus\{\alpha\},\lambda\rangle = 0$ and $\langle \alpha,\lambda\rangle \ne 0$.
We may assume $\lambda_\alpha$ is in $\Z\check{\Pi}_1$ or $\Z\check{\Pi}_2$.
Then $M'$ corresponds to $\{\alpha\in\Pi \mid \chi(\tau_{\lambda_\alpha}) = 0\}$~\cite[Proof of Corollary~1.5]{arXiv:0910.4570}.
If $\alpha\in\Pi_1$, then $\tau_{\lambda_\alpha}\in\Z\check{\Pi}_1\subset X_{H_1,*,+}$.
Therefore, there exists $\chi_1'\in\Satakepar(\pi_1|_{H_1})$ and $\chi_2'\in\Satakepar(\pi_2|_{H_1})$ such that $\chi(\tau_{\lambda_\alpha}) = \chi_1'(\tau_{\lambda_\alpha})\chi_2'(\tau_{\lambda_\alpha})$ by Proposition~\ref{prop:Satake parameters of tensor product}.
Since $\pi_2|_{H_1}$ is a direct sum of characters, $\chi_2(\tau_{\lambda_\alpha})\ne 0$ by Proposition~\ref{prop:Satake parameter of a character}.
Hence $\chi(\tau_{\lambda_\alpha}) = 0$ if and only if $\chi_1'(\tau_{\lambda_\alpha}) = 0$.
By Proposition~\ref{prop:comparison of Satake parameters, derived group}, $\Satakepar(\pi_1|_{H_1}) = \Satakepar(\pi_1|_{M_1})|_{\overline{\kappa}[X_{H_1,*,+}]} = \{\chi_1\}|_{\overline{\kappa}[X_{H_1,*,+}]}$.
Therefore, we have $\chi_1'(\tau_{\lambda_\alpha}) = \chi_1(\tau_{\lambda_\alpha})$.
It is zero if and only if $\alpha\in\Pi_{M}\cap\Pi_1$.
By the same argument, for $\alpha\in\Pi_2$, $\chi(\tau_{\lambda_\alpha}) = 0$ if and only if $\alpha\in\Pi_{M}\cap\Pi_2$.
Hence $M' = M$.
\end{proof}

\begin{cor}\label{cor:Satake parameter of tensor product, singular tensor trivial Satake par}
Assume that $M = M_1$.
Then we have $\Satakepar(\pi) = \Satakepar(\pi_1)\otimes\Satakepar(\pi_2) = \{(M_1,\chi_{M_1}\chi_T|_{X_{M_1,*,0}})\}$.
\end{cor}
\begin{proof}
Take $\chi\in\Satakepar(\pi)$ and let $\chi_M\colon X_{M,*,0}\to\overline{\kappa}^\times$ such that $\chi$ is parameterized by $(M,\chi_M)$.
The character $\chi_M$ is given by a restriction of $\chi$ on $X_{*,+} \cap \Pi_M^\perp = X_{*,+}\cap\Pi_1^\perp = X_{H_2,*,+}$.
Therefore, $\chi = \chi_1\otimes\chi_2$.
\end{proof}

\section{A theorem of changing the weight}\label{sec:A theorem of changing the weight}
In this section, we assume that the derived group of $G$ is simply connected.
For $\alpha\in\Pi$, we denote a fundamental weight corresponding to $\alpha$ by $\omega_\alpha$.
\subsection{Changing the weight}\label{subsec:Changing the weight}
We prove the following theorem which is a generalization of Herzig's theorem~\cite[Corollarry~6.11]{arXiv:1005.1713}.
\begin{thm}\label{thm:changing the weight}
Let $V_1,V_2$ be irreducible representations of $K$ with lowest weight $\nu_1,\nu_2$, respectively.
Assume that $\langle \nu_1,\check{\alpha}\rangle = 0$ and $\nu_2 = \nu_1 - (q - 1)\omega_\alpha$ for some $\alpha\in\Pi$.
Let $\chi\colon \overline{\kappa}[X_{*,+}]\to \overline{\kappa}$ be an algebra homomorphism parameterized by $(M,\chi_M)$.
Assume that $\alpha\not\in \Pi_M$.
If $\check{\alpha}\not\in X_{M,*,0}$ or $\chi_M(\check{\alpha})\ne 1$, then
\[
	\cInd_K^GV_1\otimes_{\mathcal{H}_G(V_1)}\chi
	\simeq
	\cInd_K^GV_2\otimes_{\mathcal{H}_G(V_2)}\chi.
\]
\end{thm}
Let $V_1,V_2,\nu_1,\nu_2$ be as above.
Fix $\lambda\in X_{*,+}$ such that $\langle \lambda,\Pi\setminus\{\alpha\}\rangle = 0$ and $\langle \lambda,\alpha\rangle \ne 0$.
Then there exist nonzero $\varphi_{21}\in \mathcal{H}_G(V_1,V_2)$ and $\varphi_{12}\in \mathcal{H}_G(V_2,V_1)$ whose support is $K\lambda(\varpi)K$.
By the proof of \cite[Corollary~6.11]{arXiv:1005.1713}, Theorem~\ref{thm:changing the weight} follows from the following lemma.
\begin{lem}\label{lem:calc of phi+*phi- with Satake}
We have $S_G(\varphi_{12}*\varphi_{21}) \in\overline{\kappa}^\times( \tau_{2\lambda} - \tau_{2\lambda - \check{\alpha}})$.
\end{lem}
This lemma follows from the following two lemmas by \cite[Proposition~5.1]{arXiv:1005.1713}.
\begin{lem}\label{lem:calc of phi+*phi-}
The composition $\varphi_{12}*\varphi_{21}$ is nonzero and its support is $K\lambda(\varpi)^2K$.
\end{lem}
\begin{lem}\label{lem:lem of dom elem}
For $\mu\in X_{*,+}$, if $\mu \le 2\lambda$ then $\mu = 2\lambda$ or $\mu \le 2\lambda - \check{\alpha}$.
\end{lem}

First, we prove Lemma~\ref{lem:calc of phi+*phi-}.
To prove it, we use the following lemma.
We use the argument in the proof of \cite[Proposition~6.8]{arXiv:1005.1713}.
For each $w\in W\simeq N_K(T(\mathcal{O}))/T(\mathcal{O})$, we fix a representative of $w$ and denote it by the same letter $w$.

\begin{lem}\label{lem:calc of varphi}
Let $V,V'$ be irreducible representations of $K$ with lowest weight $\nu,\nu'$, $v\in V,v'\in V'$ its lowest weight vector, respectively.
Assume that for $\mu\in X_{*,+}$, $V^{\overline{N_\mu}(\kappa)}\simeq (V')^{\overline{N_\mu}(\kappa)}$ as $M_\mu(\kappa)$-representations.
Let $\varphi\in \mathcal{H}_G(V,V')$ be such that $\supp\varphi = K\mu(\varpi)K$ and $\varphi(\mu(\varpi))v = v'$.
Put $\overline{I} = \red^{-1}(\overline{B}(\kappa))$ and $t = \mu(\varpi)$.
Then we have
\[
	\varphi*[1,v]
	=
	\sum_{w\in W_{\nu}/(W_{\nu}\cap W_\mu)}\sum_{a\in (w^{-1}\overline{I}w\cap \overline{N}(\mathcal{O}))/t^{-1}\overline{N}(\mathcal{O})t}[wat^{-1},v'].
\]
\end{lem}
\begin{proof}
We have
\[
	(\varphi*[1,v])(x) = \sum_{y\in G/K}\varphi(y)[1,v](xy) = \sum_{y\in KtK/K}\varphi(y)[1,v](xy). 
\]
If this is not zero, then $xy\in K$ for some $y\in KtK$.
Hence $x\in Kt^{-1}K$.
Namely, $\supp(\varphi*[1,v]) \subset Kt^{-1}K$.
The value at $x = kt^{-1}$ for $k\in K$ is
\[
	(\varphi*[1,v])(kt^{-1}) = \sum_{y\in KtK/K}\varphi(y)[1,v](kt^{-1}y) = \varphi(t)[1,v](k) = \varphi(t)k^{-1}v.
\]
Therefore, we have
\[
	\varphi*[1,v] = \sum_{k\in K/(K\cap t^{-1}Kt)}[kt^{-1},\varphi(t)k^{-1}v].
\]

Put $P = P_\mu$.
We have $K\cap t^{-1}Kt\supset P(\mathcal{O})$ and $\red(K\cap t^{-1}Kt) = P(\kappa)$.
Therefore, we have a surjective map $G(\mathcal{O})/P(\mathcal{O})\twoheadrightarrow K/(K\cap t^{-1}Kt)$.
For each $w\in W\simeq N_K(T(\mathcal{O}))/T(\mathcal{O})$, we fix a representative of $w$ and denote it by the same letter $w$.
Then we have
\[
	G(\mathcal{O})
	=
	\coprod_{w\in W/W_\mu}
	w(w^{-1}\overline{I}w\cap \overline{N}(\mathcal{O})) P(\mathcal{O}).
\]
Hence $\varphi*[1,v]$ is a sum of a form $[wat^{-1},\varphi(t)a^{-1}w^{-1}v]$ for $a\in w^{-1}\overline{I}w\cap\overline{N}(\mathcal{O})$ and $w\in W/W_\mu$.
We prove that $\varphi(t)a^{-1}w^{-1}v \ne 0$ implies $w\in W_{\nu}W_\mu$.
Since $\red(a)\in w^{-1}\overline{B}(\kappa)w\cap \overline{N}(\kappa)\subset  w^{-1}\overline{U}(\kappa)w$, we have $a^{-1}w^{-1}v = w^{-1}v$.
The homomorphism $\varphi(t)$ is given by $V\twoheadrightarrow (V)_{N_\mu(\kappa)}\simeq (V')^{\overline{N_\mu}(\kappa)}\hookrightarrow V'$.
Hence if $\varphi(t)w^{-1}v\ne 0$, then $w^{-1}v\in V^{\overline{N_\mu}(\kappa)}$.
Since $\{g\in G(\kappa)\mid gv\in \overline{\kappa}v\} = \overline{P_{\nu}}(\kappa)$, we have $\overline{P_{\nu}}(\kappa) \supset w\overline{N_\mu}(\kappa)w^{-1}$.
Then $\Delta_{\nu}^+\cup\Delta^+\supset w(\Delta^+\setminus\Delta_\mu^+)$.
Hence, $(\Delta^-\setminus\Delta_{\nu}^-)\cap w(\Delta^+\setminus\Delta_\mu^+) = \emptyset$.
Take $w'\in W_{\nu} w W_\mu$ such that $w'$ is shortest in $W_{\nu} w W_\mu$~\cite[Ch.~IV, Exercises, \S 1~(3)]{MR1890629}.
Then $(\Delta^-\setminus\Delta_{\nu}^-)\cap w'(\Delta^+\setminus\Delta_\mu^+) = \emptyset$.
By the condition of $w'$, $\Delta^-\cap w'(\Delta^+\setminus\Delta_\mu^+) = \Delta^-\cap w'\Delta^+$ and $(\Delta^-\setminus\Delta_{\nu}^-)\cap w'\Delta^+ = \Delta^-\cap w'\Delta^+$.
Therefore, we have $\Delta^- \cap w'\Delta^+ = \emptyset$.
Hence $w' = 1$.
We have $w\in W_{\nu} W_\mu/W_\mu = W_{\nu}/(W_{\nu}\cap W_\mu)$.
Hence we may assume $w \in W_{\nu}$.
Therefore, $\varphi(t)w^{-1}v = \varphi(t)v = v'$.
Hence,
\[
	\varphi*[1,v]
	=
	\sum_{w\in W_{\nu}/(W_{\nu}\cap W_\mu)}\sum_{a\in (w^{-1}\overline{I}w\cap \overline{N}(\mathcal{O}))/(w^{-1}\overline{I}w\cap \overline{N}(\mathcal{O})\cap t^{-1}Kt)}[wat^{-1},v'].
\]
Since $\langle \alpha,\mu\rangle < 0$ for all weight $\alpha$ of $\overline{N}$, $t = \mu(\pi)$ satisfies $t\overline{N}(\mathcal{O})t^{-1}\supset \overline{N}(\mathcal{O})$.
Hence $t\overline{N}(\mathcal{O})t^{-1}\cap K = \overline{N}(\mathcal{O})$.
Equivalently, we have $\overline{N}(\mathcal{O})\cap t^{-1}Kt = t^{-1}\overline{N}(\mathcal{O})t$.
We also have $\red(t^{-1}\overline{N}(\mathcal{O})t)$ is trivial.
Hence $t^{-1}\overline{N}(\mathcal{O})t\subset w^{-1}\overline{I}w$.
Therefore, $w^{-1}\overline{I}w\cap \overline{N}(\mathcal{O})\cap t^{-1}Kt = t^{-1}\overline{N}(\mathcal{O})t$.
Hence we have
\[
	\varphi *[1,v]
	=
	\sum_{w\in W_{\nu}/(W_{\nu}\cap W_\mu)}\sum_{a\in (w^{-1}\overline{I}w\cap \overline{N}(\mathcal{O}))/t^{-1}\overline{N}(\mathcal{O})t}[wat^{-1},v'].
\]
\end{proof}

\begin{proof}[Proof of Lemma~\ref{lem:calc of phi+*phi-}]
Put $t = \lambda(\varpi)$.
Let $v_1\in V_1,v_2\in V_2$ be lowest weight vectors.
We may assume $\varphi_{21}(t)v_1 = v_2$ and $\varphi_{12}(t)v_2 = v_1$.
By Lemma~\ref{lem:calc of varphi}, we have
\[
	\varphi_{21} *[1,v_1]
	=
	\sum_{w\in W_{\nu_1}/(W_{\nu_1}\cap W_\lambda)}\sum_{a\in (w^{-1}\overline{I}w\cap \overline{N}(\mathcal{O}))/t^{-1}\overline{N}(\mathcal{O})t}[wat^{-1},v_2].
\]
By the assumption, $W_{\nu_2}\cap W_\lambda = W_{\nu_2}$.
Hence we have
\[
	\varphi_{12}*[1,v_2]
	=
	\sum_{b\in \overline{N}(\mathcal{O})/t^{-1}\overline{N}(\mathcal{O})t}
	[bt^{-1},v_1]
\]
by Lemma~\ref{lem:calc of varphi}.
Therefore, we have
\begin{align*}
	&\varphi_{12}*\varphi_{21}*[1,v_1]\\
	&=
	\varphi_{12}*\left(
	\sum_{w\in W_{\nu_1}/(W_\lambda\cap W_{\nu_1})}\sum_{a\in(w^{-1}\overline{I}w\cap \overline{N}(\mathcal{O}))/t^{-1}\overline{N}(\mathcal{O})t}
	[wat^{-1},v_2]
	\right)\\
	& = 
	\sum_{w\in W_{\nu_1}/(W_\lambda\cap W_{\nu_1})}\sum_{a\in(w^{-1}\overline{I}w\cap \overline{N}(\mathcal{O}))/t^{-1}\overline{N}(\mathcal{O})t}
	wat^{-1}\varphi_{12}*[1,v_2]\\
	& = 
	\sum_{w\in W_{\nu_1}/(W_\lambda\cap W_{\nu_1})}\sum_{a\in(w^{-1}\overline{I}w\cap \overline{N}(\mathcal{O}))/t^{-1}\overline{N}(\mathcal{O})t}
	\sum_{b\in \overline{N}(\mathcal{O})/t^{-1}\overline{N}(\mathcal{O})t}
	[wat^{-1}bt^{-1},v_1]\\
	& = 
	\sum_{w\in W_{\nu_1}/(W_\lambda\cap W_{\nu_1})}\sum_{a\in(w^{-1}\overline{I}w\cap \overline{N}(\mathcal{O}))/t^{-1}\overline{N}(\mathcal{O})t}
	\sum_{b\in t^{-1}\overline{N}(\mathcal{O})t/t^{-2}\overline{N}(\mathcal{O})t^2}
	[wabt^{-2},v_1]\\
	& = 
	\sum_{w\in W_{\nu_1}/(W_\lambda\cap W_{\nu_1})}\sum_{c\in(w^{-1}\overline{I}w\cap \overline{N}(\mathcal{O}))/t^{-2}\overline{N}(\mathcal{O})t^2}
	[wct^{-2},v_1].
\end{align*}
Let $\varphi\in \mathcal{H}_G(V_1)$ whose support is $K\lambda(\varpi)^2K$ and $\varphi(\lambda(\varpi)^2)v_1 = v_1$.
By Lemma~\ref{lem:calc of varphi}, the right hand side of the above equation is $\varphi *[1,v_1]$. (Notice that $W_\lambda = W_{2\lambda}$.)
Since $[1,v_1]$ generates $\cInd_K^G(V_1)$, we get the lemma.
\end{proof}

Finally, we prove Lemma~\ref{lem:lem of dom elem}.
\begin{proof}[Proof of Lemma~\ref{lem:lem of dom elem}]
Assume that $\mu\le 2\lambda$ and $\mu\not\le 2\lambda - \check{\alpha}$.
Since $\mu \le 2\lambda$, there exists $n_\beta \in\Z_{\ge 0}$ such that $2\lambda - \mu = \sum_{\beta\in\Pi}n_\beta \check{\beta}$.
Then for $\gamma\in\Pi\setminus\{\alpha\}$, we have $\sum_{\beta}n_\beta \langle \gamma,\check{\beta}\rangle = \langle \gamma,2\lambda - \mu\rangle = -\langle \gamma,\mu\rangle\le 0$.
By the assumption, $n_\alpha = 0$.
Then $\sum_{\beta\ne\alpha}n_\beta\langle \gamma,\check{\beta}\rangle \le 0$.
Since the matrix $(\langle \gamma,\check{\beta}\rangle)_{\gamma,\beta\ne\alpha}$ is positive definite, we have $n_\beta = 0$ for all $\beta\in\Pi\setminus\{\alpha\}$.
Hence $\mu = 2\lambda$.
\end{proof}

\subsection{Comparison of composition factors}
Herzig proved the following proposition when $G = \GL_n$~\cite[Proposition~9.1]{arXiv:1005.1713}.
Let $\Ord_{\overline{P}}(\pi)$ be the ordinary part of $\pi$ defined by Emerton~\cite{MR2667882}.
\begin{prop}\label{prop:herzig prop 9.1}
Let $\pi$ be an admissible representation of $G$ which contains the trivial $K$-representation.
Assume that there exists $\chi\in \Satakepar(\pi,\trivrep_K)$ which is parameterized by $(T,\trivrep_{X_{T,*,0}} = \trivrep_{X_*})$.
Then $\pi$ contains the trivial representation or $\Ord_{\overline{P}}(\pi)\ne 0$ for some proper parabolic subgroup $P$.
\end{prop}
We generalize this proposition for any $G$.
He proved this proposition by a calculation in the affine Hecke algebra.
Here, we prove the proposition in a different way.
We prove the proposition from the following proposition. (In fact, we will use only the following proposition.)
When $G = \GL_2$, this proposition is proved by Barthel-Livn\'e~\cite[Theorem~20]{MR1361556}.
\begin{prop}\label{prop:changing the weight, by character}
Let $\chi\colon\overline{\kappa}[X_*]\to \overline{\kappa}$ be an algebra homomorphism and $V$ an irreducible representation of $K$.
Then $\cInd_K^G(V)\otimes_{\mathcal{H}_G(V)}\chi$ has a finite length and its composition factors depend only on $\chi$ and the $T(\kappa)$-representation $V^{\overline{U}(\kappa)}$.
\end{prop}
For a parabolic subgroup $P\subset G$, let $\Sp_P$ be the special representation~\cite{Grosse-special-rep}.
We have the following corollary.
\begin{cor}\label{cor:composition factors, from special rep}
Let $V$ be an irreducible $K$-representation such that $V^{\overline{U}(\kappa)}$ is trivial and $\chi\colon \overline{\kappa}[X_*]\to \overline{\kappa}$ an algebra homomorphism parameterized by $(T,\trivrep_{X_*})$.
Then the composition factors of $\cInd_K^G(V)\otimes_{\mathcal{H}_G(V)}\chi$ are $\{\Sp_P\mid P\subset G\}$.
\end{cor}
\begin{proof}
Let $V_1$ be the irreducible $K$-representation with lowest weight $-\sum_{\alpha\in\Pi}(q - 1)\omega_\alpha$.
Then we have $V^{\overline{U}(\kappa)}\simeq V_1^{\overline{U}(\kappa)}\simeq \trivrep_{T(\kappa)}$.
By Proposition~\ref{prop:changing the weight, by character}, we have that $\cInd_K^G(V)\otimes_{\mathcal{H}_G(V)}\chi$ and $\cInd_K^G(V_1)\otimes_{\mathcal{H}_G(V_1)}\chi$ have the same composition factors.
By Herzig's theorem~\cite[Theorem~3.1]{arXiv:1005.1713}, we have 
\[
	\cInd_K^G(V_1)\otimes_{\mathcal{H}_G(V_1)}\chi\simeq\Ind_B^G(\cInd_{T\cap K}^T(\trivrep_{T\cap K})\otimes_{\mathcal{H}_T(\trivrep_{T\cap K})}\chi) = \Ind_B^G(\trivrep_T).
\]
Hence the corollary follows from \cite[Corollary~7.3]{arXiv:1005.1713}.
\end{proof}
Using this corollary, we can prove Proposition~\ref{prop:herzig prop 9.1}.
\begin{proof}[Proof of Proposition~\ref{prop:herzig prop 9.1}]
Let $\chi\colon \overline{\kappa}[X_{*,+}]\to \overline{\kappa}$ be an algebra homomorphism parameterized by $(T,\trivrep_{X_*})$.
Then from the assumption, we have a nonzero homomorphism $\cInd_K^G(\trivrep_K)\otimes_{\mathcal{H}_G(\trivrep_K)}\chi\to \pi$.
Hence $\pi$ contains an irreducible subquotient of $\cInd_K^G(\trivrep_K)\otimes_{\mathcal{H}_G(\trivrep_K)}\chi$ as a subrepresentation.
By Corollary~\ref{cor:composition factors, from special rep}, such subquotient is $\Sp_P$ for a parabolic subgroup $P$.
If $P = G$, then $\trivrep_G = \Sp_G\subset\pi$.
If $P\ne G$, then $0\ne \Ord_{\overline{P}}(\Sp_P)\hookrightarrow \Ord_{\overline{P}}(\pi)$.
\end{proof}
\begin{rem}
If $\pi$ is irreducible, then $\pi \simeq\Sp_P$.
Since $\pi$ contains the trivial $K$-representation, $\pi$ is trivial by \cite[Proposition~7.4]{arXiv:1005.1713}.
\end{rem}

In the rest of this section, we prove Proposition~\ref{prop:changing the weight, by character}.
We use the following theorem due to Herzig~\cite[Theorem~3.1]{arXiv:1005.1713}.
\begin{thm}\label{thm:compact induction and parabolic induction}
Let $V$ be an irreducible representation of $K$ with lowest weight $\nu$, $P = MN$ a standard parabolic subgroup.
Assume that $\Stab_W(\nu)\subset W_M$.
Then we have
\[
	\cInd_K^G(V)\otimes_{\mathcal{H}_G(V)}\mathcal{H}_M(V^{\overline{N}}(\kappa))
	\simeq
	\Ind_P^G(\cInd_{M\cap K}^MV^{\overline{N}(\kappa)})
\]
as $G$-representations and $\mathcal{H}_M(V^{\overline{N}(\kappa)})$-modules.
\end{thm}
\begin{rem}
In fact, the theorem of Herzig is weaker than this theorem.
However, his proof can be applicable for this theorem.
\end{rem}
For a parabolic subgroup $P = MN$, let $V_P$ be the irreducible representation of $K$ with lowest weight $-\sum_{\alpha\in\Pi\setminus\Pi_M}(q - 1)\omega_\alpha$.
Put $\pi_P = \Ind_K^G(V_P)\otimes_{\mathcal{H}_G(V_P)}\overline{\kappa}[X_*]$.
Then we have $\pi_P\simeq \Ind_P^G(\cInd_{M\cap K}^M(\trivrep_{M\cap K})\otimes_{\mathcal{H}_M(\trivrep_{M(\kappa)})}\overline{\kappa}[X_*])$ by Theorem~\ref{thm:compact induction and parabolic induction}. (Notice that $(V_P)^{\overline{N}(\kappa)}$ is the trivial representation.)
In particular, we have $\pi_B\simeq \Ind_B^G(\overline{\kappa}[X_*])$.
Here, $T$ acts on $\overline{\kappa}[X_*]$ by $T\to T/T(\mathcal{O})\simeq X_*\to \End(\overline{\kappa}[X_*])$. (The last map is given by the multiplication.)

\begin{lem}\label{lem:homomorphism between compact induction}
For parabolic subgroups $P\subset P'$, there exist $\Phi_{P,P'}\colon\pi_{P'}\to \pi_P$ and $\Phi_{P',P}\colon\pi_P\to\pi_{P'}$ which have the following properties.
\begin{enumerate}
\item $\Phi_{P,P'}$ and $\Phi_{P',P}$ are $G$- and $\overline{\kappa}[X_*]$-equivariant.
\item $\Phi_{P,P} = \id$.
\item For $P_1\subset P_2\subset P_3$, $\Phi_{P_1,P_2}\circ\Phi_{P_2,P_3} = \Phi_{P_1,P_3}$ and $\Phi_{P_3,P_2}\circ\Phi_{P_2,P_1} = \Phi_{P_3,P_1}$.
\item For $P\subset P'$, compositions $\Phi_{P,P'}\circ\Phi_{P',P}$ and $\Phi_{P',P}\circ\Phi_{P,P'}$ are given by $\prod_{\alpha\in \Pi_{P'}\setminus\Pi_P}(\tau_{\check{\alpha}} - 1)$.
\end{enumerate}
\end{lem}
\begin{proof}
For each $\alpha\in\Pi$, fix $\lambda_\alpha\in X_{*,+}$ such that $\langle\lambda_\alpha,\Pi\setminus\{\alpha\}\rangle = 0$ and $\langle\lambda_\alpha,\alpha\rangle\ne 0$.
We also fix a lowest weight vector $v_P$ of $V_P$.

Let $P_1\subset P_2$ be parabolic subgroups such that $\#\Pi_{P_2} = \#\Pi_{P_1} + 1$ and $\Pi_{P_2} = \Pi_{P_1}\cup\{\alpha\}$.
Take $\varphi_{P_2,P_1}\in \mathcal{H}_G(V_{P_1},V_{P_2})$ and $\varphi_{P_1,P_2}\in \mathcal{H}_G(V_{P_2},V_{P_1})$ such that their support is $K\lambda_\alpha(\varpi)K$ and their value at $\lambda_\alpha(\varpi)$ send the lowest weight vector to the lowest weight vector (as in subsection~\ref{subsec:Changing the weight}).
The elements $\varphi_{P_2,P_1}$ and $\varphi_{P_1,P_2}$ give homomorphisms $\pi_{P_1}\to\pi_{P_2}$ and $\pi_{P_2}\to \pi_{P_1}$.
Let $\Phi_{P_1,P_2}$ (resp.~$\Phi_{P_2,P_1}$) be a homomorphism given by $\varphi_{P_1,P_2}$ (resp.~$-\tau_{\check{\alpha} - 2\lambda_\alpha}\varphi_{P_2,P_1}$).
By Lemma~\ref{lem:calc of phi+*phi- with Satake}, these homomorphisms satisfy the condition (4).
For general $P'\subset P$, take a chain of parabolic subgroups $P' = P_1\subset \dotsb \subset P_r = P$ such that $\#\Pi_{P_{i + 1}} = \#\Pi_{P_i} + 1$.
Define $\Phi_{P',P} = \Phi_{P_1,P_2}\circ\dotsb\circ\Phi_{P_{r - 1},P_r}$ and $\Phi_{P,P'} = \Phi_{P_r,P_{r - 1}}\circ\dotsb\circ\Phi_{P_2,P_1}$.
The by \cite[Proposition~6.3]{arXiv:1005.1713}, the condition (4) are satisfied.

It is sufficient to prove that $\Phi_{P',P}$ and $\Phi_{P,P'}$ are independent of the choice of a chain.
To prove this, we may assume that the length of the chain is $2$.
So let $P,P',P_1,P_2$ be parabolic subgroups and $\alpha,\beta\in\Pi$ such that $\alpha\ne \beta$, $\alpha,\beta\not\in\Pi_P$, $\Pi_{P_1} = \Pi_P\cup \{\alpha\}$, $\Pi_{P_2} = \Pi_P\cup\{\beta\}$ and $\Pi_{P'} = \Pi_P\cup\{\alpha,\beta\}$.
Put $t_\alpha = \lambda_\alpha(\varpi)$ and $t_\beta = \lambda_\beta(\varpi)$.
Then by Lemma~\ref{lem:calc of varphi}, we have
\begin{align*}
(\Phi_{P',P_1}\circ\Phi_{P_1,P})([1,v_P]) & = \sum_{a\in \overline{N}(\mathcal{O})/t_\alpha^{-1}\overline{N}(\mathcal{O})t_\alpha}\Phi_{P',P_1}([at_\alpha^{-1},v_{P_1}])\\
& = \sum_{a\in \overline{N}(\mathcal{O})/t_\alpha^{-1}\overline{N}(\mathcal{O})t_\alpha}\sum_{b\in \overline{N}(\mathcal{O})/t_\beta^{-1}\overline{N}(\mathcal{O})t_\beta}[at_\alpha^{-1}bt_\beta^{-1},v_{P'}]\\
& = \sum_{c\in \overline{N}(\mathcal{O})/(t_\alpha t_\beta)^{-1}\overline{N}(\mathcal{O})(t_\alpha t_\beta)}[c(t_\alpha t_\beta)^{-1},v_{P'}].
\end{align*}
Hence we have $(\Phi_{P',P_1}\circ\Phi_{P_1,P})([1,v_P]) = (\Phi_{P',P_2}\circ\Phi_{P_2,P})([1,v_P])$.
Therefore, $\Phi_{P',P_1}\circ\Phi_{P_1,P} = \Phi_{P',P_2}\circ\Phi_{P_2,P}$.

Since $\Phi_{P',P_1}\circ\Phi_{P_1,P}$ satisfies the condition (4), 
\[
(\tau_{\check{\alpha}} - 1)(\tau_{\check{\beta}} - 1)(\Phi_{P,P_2}\circ \Phi_{P_2,P'}) = (\Phi_{P,P_2}\circ\Phi_{P_2,P'})\circ(\Phi_{P',P_1}\circ\Phi_{P_1,P}\circ\Phi_{P,P_1}\circ\Phi_{P_1,P'}).
\]
By $\Phi_{P',P_1}\circ\Phi_{P_1,P} = \Phi_{P',P_2}\circ\Phi_{P_2,P}$, the right hand side is equal to
\[
	(\Phi_{P,P_2}\circ\Phi_{P_2,P'}\circ\Phi_{P',P_2}\circ\Phi_{P_2,P})\circ(\Phi_{P,P_1}\circ\Phi_{P_1,P'}).
\]
Using the condition (4) for $\Phi_{P,P_2}\circ\Phi_{P_2,P'}$, this is equal to
\[
	(\tau_{\check{\alpha}} - 1)(\tau_{\check{\beta}} - 1)(\Phi_{P,P_1}\circ\Phi_{P_1,P'}).
\]
Since $\pi_P$ is a torsion-free $\overline{\kappa}[X_*]$-module~\cite[Corollary~6.5]{arXiv:1005.1713}, we have $\Phi_{P,P_2}\circ \Phi_{P_2,P'} = \Phi_{P,P_1}\circ \Phi_{P_1,P'}$.
We get the lemma.
\end{proof}

We fix such homomorphisms.
Since $\pi_P$ is a torsion-free $\overline{\kappa}[X_*]$-module~\cite[Corollary~6.5]{arXiv:1005.1713}, the condition (4) implies $\Phi_{P,P'}$ and $\Phi_{P',P}$ are injective.

\begin{lem}
We have $\pi_P^K\simeq \overline{\kappa}[X_*]$.
\end{lem}
\begin{proof}
We have $\pi_P\simeq \Ind_{P\cap K}^K(\cInd_{M\cap K}^M(\trivrep_{M\cap K})\otimes_{\mathcal{H}_M(\trivrep_{M\cap K})}\overline{\kappa}[X_*])$ by the Iwasawa decomposition $G = KP$.
Therefore, we have
\begin{multline*}
\pi_P^K = \Hom_K(\trivrep_K,\Ind_{P\cap K}^K(\cInd_{M\cap K}^M(\trivrep_{M\cap K})\otimes_{\mathcal{H}_M(\trivrep_{M\cap K})}\overline{\kappa}[X_*]))\\
\simeq\Hom_{M\cap K}(\trivrep_{M\cap K},\cInd_{M\cap K}^M(\trivrep_{M\cap K})\otimes_{\mathcal{H}_M(\trivrep_{M\cap K})}\overline{\kappa}[X_*])\\
\simeq\End_{M}(\cInd_{M\cap K}^M(\trivrep_{M\cap K}))\otimes_{\mathcal{H}_M(\trivrep_{M\cap K})}\overline{\kappa}[X_*]
\simeq\overline{\kappa}[X_*].
\end{multline*}
\end{proof}

\begin{rem}\label{rem:isom of K-invs of principal series}
A homomorphism $\Ind_B^G(\overline{\kappa}[X_*])\ni f\mapsto f(1)\in\overline{\kappa}[X_*]$ gives an isomorphism $\pi_B^K\simeq\overline{\kappa}[X_*]$.
\end{rem}
Set $f_0 = [1,1]\otimes 1\in \cInd_K^G(\trivrep_K)\otimes_{\mathcal{H}_G(\trivrep_K)}\overline{\kappa}[X_*] = \pi_G$.
Then $\pi_G^K$ is generated by $f_0$ as a $\overline{\kappa}[X_*]$-module.
We also have that $\pi_G$ is generated by $\pi_G^K = \overline{\kappa}[X_*]f_0$ as a $G$-module.
\begin{lem}\label{lem:hom between parabolic induction}
Let $P = MN$ be a parabolic subgroup and $\sigma_1,\sigma_2$ representations of $M$.
Then we have $\Hom_M(\sigma_1,\sigma_2)\simeq\Hom_G(\Ind_P^G(\sigma_1),\Ind_P^G(\sigma_2))$.
\end{lem}
\begin{proof}
If $P = B$, this lemma is proved by Vign\'eras~\cite[Corollaire~7]{MR2399093}.
The same proof can be applicable.
\end{proof}

\begin{lem}
The element $\tau_{\check{\alpha}} - 1\in \overline{\kappa}[X_*]$ is irreducible.
\end{lem}
\begin{proof}
Take $d\in \Z_{>0}$ and $\lambda\in X_*$ such that $\langle \alpha,X_*\rangle = d\Z$ and $\langle \alpha,\lambda\rangle = d$.
Then we have $X_* = \Z\lambda\oplus \Ker\alpha$.
Let $a,b\in\overline{\kappa}[X_*]$ such that $\tau_{\check{\alpha}} - 1 = ab$.
Put $t = \tau_\lambda$.
Then we have $a = \sum_n a_nt^n$ and $b_n = \sum_n b_nt^n$ where $a_n,b_n\in \overline{\kappa}[\Ker\alpha]$.
Put $k_a = \max\{n\mid a_n\ne 0\}$, $l_a = \min\{n\mid a_n\ne 0\}$, $k_b = \max\{n\mid b_n\ne 0\}$, $l_b = \min\{n\mid b_n\ne 0\}$.
We may assume $k_a - l_a\le k_b - l_b$.
Take $c\in\Z$ and $\lambda_0\in \Ker\alpha$ such that $\check{\alpha} = c\lambda + \lambda_0$.
Then $c = 1$ or $2$ and  we have $ab = \tau_{\check{\alpha}} - 1 = t^c\tau_{\lambda_0} - 1$.
Therefore, $k_a + k_b = c$ and $a_{k_a}b_{k_b} = \tau_{\lambda_0}\in\overline{\kappa}[\Ker\alpha]^\times$.
Replacing $(a,b)$ with $(au^{-1},bu)$ for $u = t^{k_a - 1}a_{k_a}\in\overline{\kappa}[X_*]^\times$, we may assume $k_a = 1$ and $a_{k_a} = 1$.
Hence $k_b = c - 1$.
We prove $a\in\overline{\kappa}[X_*]^\times$.
If $k_a = l_a$, then $a = t\in\overline{\kappa}[X_*]^\times$.
Hence we may assume $k_a\ne l_a$.
By $ab = \tau_{\check{\alpha}} - 1 = t^c\tau_{\lambda_0} - 1$, we have $l_a + l_b = 0$.
Therefore, $(c,k_a,l_a,k_b,l_b)$ satisfies the following conditions:
\[
	\text{$c = 1$ or $2$},\ k_a = 1,\ k_b = c - 1,\ l_a < k_a,\ k_a - l_a\le k_b - l_b,\ l_a + l_b = 0.
\]
From $k_a = 1$, $k_b = c - 1$ and $k_a - l_a \le k_b - l_b$, we have $1 - l_a \le c - 1 - l_b$.
Since $l_a + l_b = 0$, we have $1 - l_a \le c - 1 + l_a$.
Therefore, $l_a \ge 1 - c/2$.
We also have $1 = k_a > l_a$.
Hence $l_a \le 0$.
From this, $0 \ge 1 - c/2$.
Hence $c = 2$.
Therefore $0\le l_a \le 1 - c/2 = 0$.
Hence $l_a = 0$ and $l_b = -l_a = 0$.
We get $(c,k_a,l_a,k_b,l_b) = (2,1,0,1,0)$.

Now we have $a = t + a_0$ and $b = b_1t + b_0$.
Since $ab = \tau_{\lambda_0}t^2 - 1$, we have
\[
	b_1 = \tau_{\lambda_0},\ a_0b_1 + b_0 = 0\ \text{and}\ a_0b_0 = -1.
\]
By the last equation, $b_0\in \overline{\kappa}[X_*]^\times$.
Hence $b_0 \in \overline{\kappa}^\times \tau_\mu$ for some $\mu\in X_*$.
We have $\tau_{\lambda_0} = b_1 = -b_0a_0^{-1} = b_0^2$.
Therefore, $\lambda_0 = 2\mu$.
Hence $\check{\alpha} = 2(\lambda + \mu)\in 2X_*$.
This is a contradiction since we assume that the derived group of $G$ is simply connected.
\end{proof}

\begin{lem}\label{lem:Phi gives pi^K=pi^K}
The image of $f_0$ under $\Phi_{B,G}$ is a basis of $\pi_B^K$.
\end{lem}
\begin{proof}
It is sufficient to prove that $\Phi_{B,G}(\pi_G^K) = \pi_B^K$.
We prove $\Phi_{B,G}(\pi_G^K) \supset\prod_{\beta\in\Pi\setminus\{\alpha\}}(\tau_{\check{\beta}} - 1)\pi_B^K$ for all $\alpha\in\Pi$.
Then for each $\alpha\in\Pi$, there exists $a_\alpha\in\overline{\kappa}[X_*]$ such that $a_\alpha\Phi_{B,G}(f_0) = \prod_{\beta\in\Pi\setminus\{\alpha\}}(\tau_{\check{\beta}} - 1)f'_0$ where $f'_0$ is a basis of $\pi_B^K$.
Since $(\tau_{\check{\alpha}} - 1)$ are distinct irreducible elements and $\overline{\kappa}[X_*]$ is UFD, we have $\Phi_{B,G}(f_0) \in\overline{\kappa}[X_*]^\times f'_0$.
Hence the lemma is proved.

So it is sufficient to prove $\Phi_{B,G}(\pi_G^K) \supset\prod_{\beta\in\Pi\setminus\{\alpha\}}(\tau_{\check{\beta}} - 1)\pi_B^K$ for all $\alpha\in\Pi$.
Fix $\alpha\in\Pi$ and let $P$ be the parabolic subgroup corresponding to $\{\alpha\}$.
Since $\Phi_{P,G}(\pi_G^K)\supset\Phi_{P,G}(\Phi_{G,P}(\pi_P^K)) = \prod_{\beta\in\Pi\setminus\{\alpha\}}(\tau_{\check{\beta}} - 1)\pi_P^K$, it is sufficient to prove $\Phi_{B,P}(\pi_P^K) = \pi_B^K$.
By Lemma~\ref{lem:hom between parabolic induction}, $\Phi_{B,P}$ is given by a certain homomorphism $\Phi\colon \cInd_{M\cap K}^M(\trivrep_{M\cap K})\otimes_{\mathcal{H}_M(\trivrep_{M\cap K})}\overline{\kappa}[X_*]\to \Ind_{M\cap B}^M(\trivrep_{M\cap B})$.
We also have $\Phi_{P,B}$ is induced by some $\Phi'\colon \Ind_{M\cap B}^M(\trivrep_{M\cap B})\to \cInd_{M\cap K}^M(\trivrep_{M\cap K})\otimes_{\mathcal{H}_M(\trivrep_{M\cap K})}\overline{\kappa}[X_*]$.
We have $\Phi'\circ \Phi = (\tau_{\check{\alpha}} - 1)$ and $\Phi\circ\Phi' = (\tau_{\check{\alpha}} - 1)$.
Therefore, it is sufficient to prove the lemma when the semisimple rank of $G$ is one.

Now we assume that the semisimple rank of $G$ is one.
Let $\Pi = \{\alpha\}$.
Take $a,b\in\overline{\kappa}[X_*]$ such that $\Phi_{B,G}(\pi_G^K) = a\pi_B^K$, $\Phi_{G,B}(\pi_B^K) = b\pi_G^K$ and $ab = \tau_{\check{\alpha}} - 1$.
Assume $\Phi_{B,G}(\pi_G^K) \ne \pi_B^K$.
It is equivalent to $a\not\in\overline{\kappa}[X_*]^\times$.
By the above lemma, $b\in\overline{\kappa}[X_*]^\times$.
Hence $\Phi_{G,B}(\pi_B^K) = \pi_G^K$.
Since $\pi_G$ is generated by $\pi_G^K$, $\Phi_{G,B}$ is surjective.
Therefore, $\Phi_{G,B}$ is isomorphic.
Let $\chi\colon\overline{\kappa}[X_*]\to\overline{\kappa}$ be a homomorphism defined by $\chi(\tau_\lambda) = 1$ for all $\lambda\in X_*$.
Then we have $\pi_B\otimes_{\overline{\kappa}[X_*]}\chi = \Ind_B^G(\trivrep_T)$.
Hence we have $\Ind_B^G(\trivrep_T)\simeq \pi_G\otimes_{\overline{\kappa}[X_*]}\chi$.
Consider a homomorphism $\cInd_K^G(\trivrep_K)\to\trivrep_G$ defined by $f\mapsto \sum_{g\in G/K}f(g)$.
This gives a homomorphism $\pi_G\otimes_{\overline{\kappa}[X_*]}\chi\to \trivrep_G$ and induced homomorphism $(\pi_G\otimes_{\overline{\kappa}[X_*]}\chi)^K\to (\trivrep_G)^K = \trivrep_G$ is an isomorphism.
We have $\trivrep_G = (\Ind_B^G(\trivrep_T))^K \hookrightarrow\Ind_B^G(\trivrep_T)\simeq \pi_G\otimes_{\overline{\kappa}[X_*]}\chi\to \trivrep_G$.
The composition is isomorphic.
Hence $\trivrep_G$ is a direct summand of $\Ind_B^G(\trivrep_T)$.
Therefore, $\End_G(\Ind_B^G(\trivrep_T))$ has a non trivial idempotent.
However, by Lemma~\ref{lem:comparison of satake parameters, restriction to Levi}, $\End_G(\Ind_B^G(\trivrep_T))\simeq\End_T(\trivrep_T)\simeq\overline{\kappa}$.
This is a contradiction.
\end{proof}
By this lemma, $\Imm\Phi_{B,G}$ is a subrepresentation of $\pi_B$ generated by $\pi_B^K$.
For each $w\in W\simeq N_K(T(\mathcal{O}))/T(\mathcal{O})$, we fix a representative of $w$ and denote it by the same letter $w$.
For a subset $A\subset W$ of $W$, let $X_{G,A}\subset\pi_B = \Ind_B^G\overline{\kappa}[X_*]$ be a $B$-stable subspace defined by $X_{G,A} = \{f\in\pi_B\mid \supp f\subset \bigcup_{w'\in A}Bw'B/B\}$.
For $w\in W$, put $X_{G,>w} = X_{G,\{w'\in W\mid w' > w\}}$ and $X_{G,\ge w} = X_{G,\{w'\in W\mid w' \ge w\}}$.
Set $X_A = X_{G,A}$, $X_{\ge w} = X_{G,\ge w}$ and $X_{>w} = X_{G,>w}$ for $A\subset W$, $w\in W$.
Set $Y = \Phi_{B,G}(\pi_G)$, $Y_A = Y\cap X_{A}$, $Y_{>w} = Y\cap X_{>w}$ and $Y_{\ge w} = Y\cap X_{\ge w}$.
For a parabolic subgroup $P = MN$, put $W(M) = \{w\in W\mid w(\Pi_M)\subset\Delta^+\}$.
Then $W(M)\times W_M\to W$ is bijective~\cite[Ch.~IV, Exercises, \S 1~(3)]{MR1890629}.

\begin{lem}
Let $P = MN$ be a parabolic subgroup, $w,v_0\in W(M)$ and $v_1\in W_M$.
Then $v_0v_1\ge w$ if and only if $v_0\ge w$.
\end{lem}
\begin{proof}
Put $v = v_0v_1$.
Let $\ell$ be the length function of $W$.
Then $\ell(v) = \ell(v_0) + \ell(v_1)$~\cite[Ch.~IV, Exercises, \S 1~(3)]{MR1890629}.
Hence $v\ge v_0$.
Therefore, $v_0\ge w$ implies $v\ge w$.

We prove $v\ge w$ implies $v_0\ge w$ by induction on $\ell(v_1)$.
If $\ell(v_1) = 0$, then $v_1 = 1$.
Hence there is nothing to prove.
Assume that $\ell(v_1) > 0$ and take $\alpha\in \Pi_M$ such that $v_1s_\alpha < v_1$ where $s_\alpha\in W_M$ is the reflection corresponding to $\alpha$.
Put $s = s_\alpha$.
Then $\ell(v_0v_1s) = \ell(v_0) + \ell(v_1s) = \ell(v_0) + \ell(v_1) - 1 = \ell(v_0v_1) - 1$.
Hence $vs < v$.
By the definition of $W(M)$, we have $ws > w$.
Hence we get $vs\ge w$~\cite[Theorem~1.1 (II, ii)]{MR0435249}.
Therefore, $v_0(v_1s) \ge w$.
Since $\ell(v_1s) < \ell(v_1)$, we have $v_0\ge w$ by inductive hypothesis.
\end{proof}

\begin{lem}\label{lem:Bruhat filtration of cInd_K^G}
We have $Y_{\ge w}/Y_{>w} = \prod_{\alpha\in\Pi,ws_\alpha < w}(\tau_{\check{\alpha}} - 1)(X_{\ge w}/X_{>w})$.
\end{lem}
\begin{proof}
Set $\Theta = \{\alpha\in\Pi\mid ws_\alpha < w\}$ and put $I = \prod_{\alpha\in\Theta}(\tau_{\check{\alpha}} - 1)\overline{\kappa}[X_*]$.
First we prove $Y_{\ge w}/Y_{>w} \subset I(X_{\ge w}/X_{>w})$, namely, we prove $Y_{\ge w}\subset IX_{\ge w} + X_{>w}$.
If $\Theta = \emptyset$, then there is nothing to prove.
So we may assume $\Theta\ne\emptyset$.
Let $P = MN$ be a parabolic subgroup corresponding to $\Theta$.
Recall that $T$ acts on $\overline{\kappa}[X_*]$ and $\pi_B = \Ind_B^G(\overline{\kappa}[X_*])$.
This action induces the action of $T$ on $\overline{\kappa}[X_*]/I$.
For $\alpha\in \Theta$, $\Imm\check{\alpha}$ acts on $\overline{\kappa}[X_*]/I$ trivially.
Therefore, the action of $T$ on $\overline{\kappa}[X_*]/I$ is extended to the action of $M$ such that $[M,M]$ acts on it trivially by Lemma~\ref{lem:G/[M_2,M_2]simeq M_1/L_2}.
We have $\Ind_P^G(\overline{\kappa}[X_*]/I)\subset\Ind_B^G(\overline{\kappa}[X_*]/I) = \pi_B/I\pi_B$.

Let $f\in (\pi_B/I\pi_B)^K = (\Ind_B^G(\overline{\kappa}[X_*]/I))^K$.
We prove $f\in \Ind_P^G(\overline{\kappa}[X_*]/I)$, namely, $f(gp) = p^{-1}f(g)$ for $g\in G$ and $p\in P$.
By the Iwasawa decomposition $G = KP$, there exist $k\in K$ and $p'\in P$ such that $g = kp'$.
Since $P = MN = [M,M]TN = ([M,M]\cap K)([M,M]\cap B)TN = ([M,M]\cap K)B$, there exist $k'\in [M,M]\cap K$ and $b\in B$ such that $p'p = k'b$.
Hence $f(gp) = f(kp'p) = f(kk'b) = b^{-1}f(1)$.
Since $k'\in [M,M]$, we have $(k')^{-1}f(1) = f(1)$.
Hence $f(gp) = (k'b)^{-1}f(1) = (p'p)^{-1}f(1)$.
Applying this to $g = 1$ and $p = p'$, we have $f(p') = (p')^{-1}f(1)$.
Hence $f(gp) = p^{-1}f(p') = p^{-1}f(kp') = p^{-1}f(g)$.
So $f\in \Ind_P^G(\overline{\kappa}[X_*]/I)$.
Hence the image of $\Phi_{G,B}(f_0)$ is in $\Ind_P^G(\overline{\kappa}[X_*]/I)$.
Therefore, the image of $Y$ is contained in $\Ind_P^G(\overline{\kappa}[X_*]/I)$.

For $f\in \pi_B$, let $\overline{f}$ be the image of $f$ under the canonical projection $\pi_B\to \pi_B/I\pi_B = \Ind_B^G(\overline{\kappa}[X_*]/I)$.
Let $f\in Y_{\ge w}$.
Then $\supp\overline{f}\subset\bigcup_{w'\ge w}Bw'B/B$.
Since $\overline{f}\in\Ind_P^G(\overline{\kappa}[X_*]/I)$, its support is right $P$-invariant.
Hence if $\supp\overline{f}\cap BwB/B\ne \emptyset$, $\supp\overline{f}\cap Bww'B/B\ne \emptyset$ for all $w'\in W_M$.
Hence $\supp\overline{f}\cap Bws_\alpha B/B \ne\emptyset$ for all $\alpha\in\Theta$.
By the definition of $\Theta$, each $\alpha\in\Theta$ satisfies $ws_\alpha < w$.
This contradicts to $\supp\overline{f}\subset\bigcup_{w'\ge w}Bw'B/B$.
So we have $\supp\overline{f}\subset\bigcup_{w'> w}Bw'B/B$.
Hence $f\in X_{>w} + I\pi_B$.

We prove $Y_{\ge w}/Y_{>w} \supset I(X_{\ge w}/X_{>w})$.
Let $P' = M'N'$ be a parabolic subgroup corresponding to $\Pi\setminus\Theta$.
First we prove that $\Phi_{B,P'}(\pi_{P'})\cap X_{\ge w}\to X_{\ge w}/X_{>w}$ is surjective.
For each parabolic subgroup $P_1 = M_1N_1\subset P'$, put $\pi_{M',P_1} = \Ind_{M'\cap P_1}^{M'}(\cInd_{M_1\cap K}^{M_1}\trivrep_{M_1\cap K}\otimes_{\mathcal{H}_{M_1}(\trivrep_{M_1\cap K})}\overline{\kappa}[X_*])$.
Then $\pi_{P_1} = \Ind_{P'}^G(\pi_{M',P_1})$.
By Lemma~\ref{lem:hom between parabolic induction}, for each $P_1\subset P_2\subset P'$, $\Phi_{P_1,P_2}$ and $\Phi_{P_2,P_1}$ are induced by some $\Phi^{M'}_{P_1,P_2}\colon \pi_{M',P_2}\to \pi_{M',P_1}$ and $\Phi^{M'}_{P_2,P_1}\colon \pi_{M',P_1}\to \pi_{M',P_2}$.
Such homomorphisms satisfy the conditions of Lemma~\ref{lem:homomorphism between compact induction}.
Therefore, $\Phi^{M'}_{P_1,P_2}$ induces a bijection $\pi_{M',P_2}^{M'\cap K}\simeq \pi_{M',P_1}^{M'\cap K}$ by Lemma~\ref{lem:Phi gives pi^K=pi^K}.
Put $\Phi = \Phi^{M'}_{B,P'}$.
Then $\Phi_{B,P'}(\pi_{P'}) = \Ind_{P'}^G(\Phi(\pi_{M',P'}))$.

Let $f\in \Phi_{B,P'}(\pi_{P'})$.
By the definition of $X_{\ge w}$, $f\in X_{\ge w}$ if and only if $\supp f\subset \bigcup_{v\ge w}BvB$.
For $v\in W$, take $v_0\in W(M')$ and $v_1\in W_{M'}$ such that $v = v_0v_1$.
Since $w\in W(M')$, $v\ge w$ if and only if $v_0\ge w$ by the above lemma.
Hence $\bigcup_{v\ge w}BvB = \bigcup_{v\ge w,v\in W(M')}BvW_{M'}B = \bigcup_{v\ge w,v\in W(M')}BvP'$.
Therefore, $\Phi_{B,P'}(\pi_{P'})\cap X_{\ge w} = \{f\in \Ind_{P'}^G(\Phi(\pi_{M',P'}))\mid \supp f \subset \bigcup_{v\ge w,v\in W(M')}BvP'/P'\}$.
Let $Z_{\ge w}$ be this space.
Put $Z_{>w} = \{f\in \Ind_{P'}^G(\Phi(\pi_{M',P'}))\mid \supp f \subset \bigcup_{v > w,v\in W(M')}BvP'/P'\}$.
Then $\Phi_{B,P'}(\pi_{M',P'})\cap X_{\ge w}\to X_{\ge w}/X_{>w}$ induces $Z_{\ge w}/Z_{>w} \to X_{\ge w}/X_{>w}$.
By the Bruhat decomposition $G/P' = \bigcup_{v\in W(M')}BvP'/P'$, the space $Z_{\ge w}/Z_{>w}$ is isomorphic to the space of locally constant compact support $\Phi(\pi_{M',P'})$-valued functions on $BwP'/P'\simeq BwB/B$.
The space $X_{\ge w}/X_{>w}$ is isomorphic to the space of locally constant compact support $\overline{\kappa}[X_*]$-valued functions on $BwB/B$.
The homomorphism $Z_{\ge w}/Z_{>w} \to X_{\ge w}/X_{>w}$ is induced by $\Phi(\pi_{M',P'})\hookrightarrow \pi_{M',B} \to \pi_{M',B}/X_{M',>1}\simeq \overline{\kappa}[X_*]$.
By Remark~\ref{rem:isom of K-invs of principal series}, $\pi_{M',B'}^{M' \cap K}\hookrightarrow \pi_{M',B'} \to \pi_{M',P'}/X_{M',>1}\simeq \overline{\kappa}[X_*]$ is isomorphic.
Hence $\Phi(\pi_{M',P'})\hookrightarrow \pi_{M',B'} \to \pi_{M',P'}/X_{M',>1}\simeq \overline{\kappa}[X_*]$ is surjective by Lemma~\ref{lem:Phi gives pi^K=pi^K}.
Therefore $\Phi_{B,P'}(\pi_{P'})\cap X_{\ge w}\to X_{\ge w}/X_{>w}$ is surjective.

We get $(\Phi_{B,P'}(\pi_{P'})\cap X_{\ge w}) + X_{>w} = X_{\ge w}$.
Since $I\Phi_{B,P'}(\pi_{P'}) = \Phi_{B,P'}(I\pi_{P'})= \Phi_{B,P'}(\Phi_{P',G}(\Phi_{G,P'}(\pi_{P'}))) = \Phi_{B,G}(\Phi_{G,P'}(\pi_{P'}))\subset \Phi_{B,G}(\pi_G) = Y$, $IX_{\ge w}\subset Y\cap X_{\ge w} + IX_{>w}\subset Y_{\ge w} + X_{>w}$.
We get the lemma.
\end{proof}
From this lemma, we get the following proposition.
\begin{prop}\label{prop:freeness of compact induction (with localization)}
Let $V$ be an irreducible representation of $K$.
The module $\cInd_K^G(V)\otimes_{\mathcal{H}_G(V)}\overline{\kappa}[X_*]$ is free as a $\overline{\kappa}[X_*]$-module.
\end{prop}
When $G = \GL_2$, this proposition is proved by Barthel-Livn\'e.
In fact, they proved that $\cInd_K^G(V)$ is a free $\mathcal{H}_G(V)$-module~\cite[Theorem~19]{MR1290194}.

\begin{proof}
Let $\nu$ be a lowest weight of $V$.
By Theorem~\ref{thm:compact induction and parabolic induction}, $\cInd_K^G(V)\otimes_{\mathcal{H}_G(V)}\overline{\kappa}[X_*]\simeq \Ind_{P_\nu}^G(\cInd_{M_\nu\cap K}^M(V^{\overline{N_\nu}(\kappa)})\otimes_{\mathcal{H}_{M_\nu}(V^{\overline{N_\nu}(\kappa)})}\overline{\kappa}[X_*])$.
It is sufficient to prove $\cInd_{M_\nu\cap K}^M(V^{\overline{N_\nu}(\kappa)})\otimes_{\mathcal{H}_{M_\nu}(V^{\overline{N_\nu}(\kappa)})}\overline{\kappa}[X_*]$ is free.
Hence we may assume $P_\nu = G$.
Therefore, $V$ is a character of $K$.
By Corollary~\ref{cor:char of K and G}, there exists a character $\nu_G$ of $G$ such that $\nu_G|_K\simeq V$.
Then $\varphi\mapsto \varphi_{\nu_G^{-1}}$ gives an isomorphism $\mathcal{H}_G(V)\simeq\mathcal{H}_G(\trivrep_K)$ (see \ref{subsec:stake parameters:definition and some lemmas}).
By this isomorphism, we can identify $\mathcal{H}_G(V)$ and $\mathcal{H}_G(\trivrep_K)$.
Under this identification, we have $\cInd_K^G(V)\otimes \nu_G^{-1}\simeq \cInd_K^G(\trivrep_K)$.
Hence we may assume $V = \trivrep_K$.
Therefore, $\cInd_K^G(V)\otimes_{\mathcal{H}_G(V)}\overline{\kappa}[X_*] = \pi_G\simeq Y$.
Since $X_{\ge w}/X_{>w}$ is free~~\cite[Lemma~3]{MR2399093}, $Y_{\ge w}/Y_{>w}$ is free by Lemma~\ref{lem:Bruhat filtration of cInd_K^G}.
Hence $Y$ is free.
\end{proof}

\begin{proof}[Proof of Proposition~\ref{prop:changing the weight, by character}]
We prove the proposition by induction on $\#\Pi_\nu$.
Namely, we prove the following by induction on $n$: If $\nu$ satisfies $\#\Pi_\nu \le n$ then the module $\cInd_K^G(V)\otimes_{\mathcal{H}_G(V)}\chi$ has a finite length and its composition factors depend only on $\chi$ and the $T(\kappa)$-representation $V^{\overline{U}(\kappa)}$.

If $\Pi_\nu = \emptyset$, then $\cInd_K^G(V)\otimes_{\mathcal{H}_G(V)}\chi$ is isomorphic to a principal series representation~\cite[Theorem~3.1]{arXiv:1005.1713}.
Hence the proposition follows.

Assume $\Pi_\nu\ne\emptyset$ and take $\alpha\in\Pi_\nu$.
Put $\nu' = \nu - (q - 1)\omega_\alpha$ and let $V'$ be the irreducible $K$-representation with lowest weight $\nu'$.
By inductive hypothesis, $\cInd_K^G(V')\otimes_{\mathcal{H}_G(V')}\chi$ has a finite length.
Define $\chi'\colon\overline{\kappa}[X_*]\to \overline{\kappa}[t,t^{-1}]$ by $\chi'(\tau_\lambda) = \chi(\tau_\lambda)t^{\langle\omega_\alpha,\lambda\rangle}$ for $\lambda\in X_*$.
(Here, $t$ is an indeterminant.)
Then $\chi$ factors through $\chi'$.
Put $\pi = \cInd_K^G(V)\otimes_{\mathcal{H}_G(V)}\chi'$ and $\pi' = \cInd_K^G(V')\otimes_{\mathcal{H}_G(V')}\chi'$.
These are free $\overline{\kappa}[t,t^{-1}]$-modules by Proposition~\ref{prop:freeness of compact induction (with localization)}.
Take $\lambda\in X_*$ such that $\langle \lambda,\Pi\setminus\{\alpha\}\rangle = 0$ and $\langle \lambda,\alpha\rangle \ne 0$.
Put $a = \chi(\tau_{\check{\alpha}})$.
As in \ref{subsec:Changing the weight}, $\lambda$ gives $\Phi\colon\pi\to \pi'$ and $\Phi'\colon\pi'\to \pi$ such that $\Phi\circ\Phi' = (at - 1)$.
Therefore, $\Phi'$ is injective and $\Imm\Phi'\subset (at - 1)\pi$.
By \cite[Lemma~2.3.4]{MR1433132}, $\pi/(t - 1)\pi$ has a finite length and $\pi/(t - 1)\pi$ and $\pi'/(t - 1)\pi'$ have the same composition factors.
\end{proof}

\section{Classification Theorem}\label{sec:Classification Theorem}
Using results in Section~\ref{sec:Satake parameters} and Section~\ref{sec:A theorem of changing the weight}, we prove the main theorem.
Almost all the proof of the theorem is a copy of Herzig's proof.
\subsection{Construction of representations}
We recall the definition of supersingular representations.
\begin{defn}[Herzig~{\cite[Definition~4.7]{arXiv:1005.1713}}]\label{defn:supersingular representation}
Let $\pi$ be an irreducible admissible representation of $G$.
\begin{enumerate}
\item The representation $\pi$ is called \emph{supersingular with respect to $(K,T,B)$} if each $\chi\in \Satakepar(\pi)$ corresponds to $(G,\chi_G)$ for some $\chi_G \colon X_{G,*,0} \to \overline{\kappa}^\times$.
\item The representation $\pi$ is called \emph{supersingular} if it is supersingular with respect to all $(K,T,B)$.
\end{enumerate}
\end{defn}
It will be proved that $\pi$ is supersingular if and only if $\pi$ is supersingular with respect to $(K,T,B)$ for a fixed $(K,T,B)$.

Now we introduce the set of parameters $\mathcal{P} = \mathcal{P}_G$.
It will parameterize the isomorphism classes of irreducible admissible representations.
Before to give $\mathcal{P}$, we give one notation.
Let $M$ be the Levi subgroup of a standard parabolic subgroup and $\sigma$ its representation with the central character $\omega_\sigma$.
Then set $\Pi_{\sigma} = \{\alpha\in\Pi\mid\langle \Pi_M,\check{\alpha}\rangle = 0,\ \omega_\sigma\circ\check{\alpha} = \trivrep_{\GL_1(F)}\}$.

Let $\mathcal{P} = \mathcal{P}_G$ be the set of $\Lambda = (\Pi_1,\Pi_2,\sigma_1)$'s such that:
\begin{itemize}
\item $\Pi_1$ and $\Pi_2$ are subsets of $\Pi$.
\item $\sigma_1$ is an irreducible admissible representation of $M_{\Pi_1}$ with the central character $\omega_{\sigma_1}$ which is supersingular with respect to $(M_{\Pi_1}\cap K,T,M_{\Pi_1}\cap B)$.
\item $\Pi_2\subset\Pi_{\sigma_1}$.
\end{itemize}
For $\Lambda = (\Pi_1,\Pi_2,\sigma_1)\in\mathcal{P}$, we attach the representation $I(\Lambda)$ of $G$ by the following way.
Let $P_\Lambda = M_\Lambda N_\Lambda$ be the standard parabolic subgroup corresponding to $\Pi_1\cup\Pi_{\sigma_1}$.
By Lemma~\ref{lem:G/[M_2,M_2]simeq M_1/L_2}, there exists the unique extension of $\sigma_1$ to $M_\Lambda$ such that $[M_{\Pi_{\sigma_1}},M_{\Pi_{\sigma_1}}]$ acts on it trivially.
We denote this representation by the same letter $\sigma_1$.
By the definition, $\Pi_1\cup\Pi_2$ is a subset of $\Pi_1\cup \Pi_{\sigma_1}$.
Hence this set defines a standard parabolic subgroup of $M_\Lambda$.
Let $\sigma_{\Lambda,2}$ be the special representation of $M_\Lambda$ with respect to this parabolic subgroup.
By the definition of special representations, $[M_{\Pi_1},M_{\Pi_1}]$ acts on $\sigma_{\Lambda,2}$ trivially and the restriction of $\sigma_{\Lambda,2}$ to $[M_{\Pi_{\sigma_1}},M_{\Pi_{\sigma_1}}]$ is the special representation of $[M_{\Pi_{\sigma_1}},M_{\Pi_{\sigma_1}}]$ with respect to the standard parabolic subgroup corresponding to $\Pi_2$.
In particular, the restriction of $\sigma_{\Lambda,2}$ to $[M_{\Pi_{\sigma_1}},M_{\Pi_{\sigma_1}}]$ is irreducible.
Put $\sigma_\Lambda = \sigma_1\otimes \sigma_{\Lambda,2}$ and $I(\Lambda) = I_G(\Lambda) = \Ind_{P_\Lambda}^G(\sigma_\Lambda)$.
By the following lemma, $\sigma_\Lambda$ is irreducible.
(Apply for $H = M_\Lambda$ and  $H' = [M_{\Pi_{\sigma_1}},M_{\Pi_{\sigma_1}}]$.)
\begin{lem}\label{lem:general lemma, tensor is irreducible}
Let $H$ be a group, $H'$ a normal subgroup of $H$ and $\sigma_2$ a representation of $H$ which is irreducible as a representation of $H'$ and $\End_{H'}(\sigma_2) = \overline{\kappa}$.
For a representation $\sigma$ of $H$, $\Hom_{H'}(\sigma_2,\sigma)$ has a structure of a representation of $H/H'$ defined by $(h\psi)(v) = h\psi(h^{-1}v)$ for $h\in H$, $\psi\in \Hom_{H'}(\sigma_2,\sigma)$ and $v\in \sigma_2$.
\begin{enumerate}
\item The natural homomorphism $\Hom_{H'}(\sigma_2,\sigma)\otimes\sigma_2\to \sigma$ is injective.
\item If $\sigma$ is irreducible, then $\Hom_{H'}(\sigma_2,\sigma)$ is zero or irreducible.
\item For an irreducible representation $\sigma_1$ of $H/H'$, $\sigma_1\otimes\sigma_2$ is an irreducible $H$-representation.
\end{enumerate}
\end{lem}

\begin{proof}
(1)
Assume that the kernel of the homomorphism is non-zero.
Take a finite dimensional subspace $V\subset\Hom_{H'}(\sigma_2,\sigma)$ such that $V\otimes\sigma_2\to \sigma$ is not injective.
This is a $H'$-homomorphism.
Therefore, there exists a subspace $V_1$ of $V$ such that the kernel is $V_1\otimes\sigma_2$.
This means $V_1 = 0$ in $\Hom_{H'}(\sigma_2,\sigma)$.
This is a contradiction.

(2)
Assume that $\sigma$ is irreducible and $\Hom_{H'}(\sigma_2,\sigma)\ne 0$.
Then by (1), we have an injective homomorphism $\Hom_{H'}(\sigma_2,\sigma)\otimes\sigma_2\hookrightarrow\sigma$.
Since $\sigma$ is irreducible, we have $\Hom_{H'}(\sigma_2,\sigma)\otimes\sigma_2\simeq\sigma$.
Therefore, $\Hom_{H'}(\sigma_2,\sigma)$ is irreducible.

(3)
Let $\sigma\subset\sigma_1\otimes\sigma_2$ be a nonzero subrepresentation.
As a representation of $H'$, $\sigma_1\otimes\sigma_2$ is a direct sum of $\sigma_2$.
Hence $\Hom_{H'}(\sigma_2,\sigma)\ne 0$.
Since $\End_{H'}(\sigma_2) = \overline{\kappa}$, we have $\Hom_{H'}(\sigma_2,\sigma_1\otimes\sigma_2) \simeq\sigma_1$.
This is an isomorphism between $H/H'$-representations.
Therefore, $\Hom_{H'}(\sigma_2,\sigma)\subset\sigma_1$.
Since $\sigma_1$ is irreducible, we have $\Hom_{H'}(\sigma_2,\sigma) = \sigma_1$.
Therefore, $\sigma = \sigma_1\otimes\sigma_2$.
\end{proof}

We use the following lemma.
It follows from Proposition~\ref{prop:Satake parameters of parbolic induction} and Corollary~\ref{cor:Satake parameter of tensor product, singular tensor trivial Satake par}.
\begin{lem}\label{lem:Satake parameter of I(lambda)}
We have $\Satakepar(I(\Lambda)) = \{(M_{\Pi_1},\chi_{\omega_{\sigma_1}})\}$, here, $\chi_{\omega_{\sigma_1}}\colon X_{M_{\Pi_1},*,0}\to\overline{\kappa}^\times$ is defined by $\chi_{\omega_{\sigma_1}}(\lambda) = \omega_{\sigma_1}(\lambda(\varpi))$.
\end{lem}

\subsection{Irreducibility of the representation}
In this subsection, we assume that the derived group of $G$ is simply connected.
We prove the irreducibility of $I(\Lambda)$.
We need a lemma.
\begin{lem}\label{lem:Stake parameter of I(Lambda), K-part}
Let $\Lambda = (\Pi_1,\Pi_2,\sigma_1)\in\mathcal{P}$, $V$ an irreducible representation of $K$ and $\nu$ its lowest weight.
Assume that $\Hom_{K}(V,I(\Lambda))\ne 0$ and $\alpha\in\Pi$ satisfies $\langle\Pi_{M_{\Pi_1}},\check{\alpha}\rangle = 0$.
Then we have $\omega_{\sigma_1}\circ \check{\alpha}|_{\mathcal{O}^\times} = \nu\circ\check{\alpha}$.
\end{lem}
\begin{proof}
Set $V_1 = V^{\overline{N_\Lambda}(\kappa)}$.
Then $V_1$ is an irreducible representation of $M_\Lambda\cap K$ with lowest weight $\nu$.
Moreover, we have $\Hom_{M_\Lambda\cap K}(V_1,\sigma_\Lambda) \ne 0$.

Let $Q$ be the parabolic subgroup of $M_\Lambda$ corresponding to $\Pi_1\cup\Pi_2$.
Then we have $\sigma_{\Lambda,2} = \Sp_{Q,M_\Lambda}$.
Put $L = [M_{\Pi_{\sigma_1}},M_{\Pi_{\sigma_1}}]$.
Then $\sigma_{\Lambda,2}|_L = \Sp_{Q\cap L,L}$.
Put $\sigma_2 = \sigma_{\Lambda,2}$ and $M_1 = M_{\Pi_1}$.

Fix $\psi\in \Hom_{M_\Lambda\cap K}(V_1,\sigma_\Lambda)\setminus\{0\}$ and consider $V_1$ as a subspace of $\sigma_\Lambda$.
Let $v\in V_1$ be a lowest weight vector.
Then we have $v\in \sigma_\Lambda^{\overline{I}_{M_\Lambda,1}}$ where $\overline{I}_{M_\Lambda,1}$ is the inverse image of $(M_\Lambda\cap \overline{U})(\kappa)$ in $M_\Lambda\cap K$.
Since $L$ acts on $\sigma_1$ trivially, we have $v\in \sigma_\Lambda^{\overline{I}_{M_\Lambda,1}}\subset \sigma_\Lambda^{\overline{I}_{M_\Lambda,1}\cap L} = \sigma_1\otimes \sigma_2^{\overline{I}_{M_\Lambda,1}\cap L}$.
Let $\overline{\sigma_2}$ be the special representation of $M_\Lambda(\kappa)$ with respect to the parabolic subgroup $Q(\kappa)$.
Then we have $\overline{\sigma_2}\hookrightarrow \sigma_2$ and we have $\overline{\sigma_2}^{(\overline{U}\cap L)(\kappa)} = \sigma_2^{\overline{I}_{M_\Lambda,1}\cap L}$~\cite{Grosse-special-rep}.
Since $\langle\Pi_{\sigma_1},\check{\Theta}_1\rangle = 0$, we have $\overline{U}\cap M_\Lambda \simeq (\overline{U}\cap L)\times (U\cap [M_1,M_1])$ as algebraic groups.
By the construction, $[M_1,M_1](\kappa)$ acts on $\overline{\sigma_2}$ trivially.
Hence we have $\overline{\sigma_2}^{(\overline{U}\cap L)(\kappa)} = \overline{\sigma_2}^{(\overline{U}\cap M_\Lambda)(\kappa)}$.
A calculation of Gro{\ss}e-K\"onne~\cite{Grosse-special-rep} shows that $T(\kappa)$ acts on $\overline{\sigma_2}^{(\overline{U}\cap M_\Lambda)(\kappa)}$ trivially.
Hence $T(\mathcal{O})$ acts on $\sigma_2^{\overline{I}_{M_\Lambda,1}\cap L}$ trivially.

Take $\alpha$ as in the lemma.
Then $\Imm\check{\alpha}\subset Z_{M_1}$.
Hence for $t\in \mathcal{O}^\times$, $\check{\alpha}(t)$ acts on $\sigma_1$ by the scalar $\omega_{\sigma_1}(\check{\alpha}(t))$.
By the above argument, $\check{\alpha}(t)$ acts on $\sigma_2^{\overline{I}_{M_\Lambda,1}\cap L}$ trivially.
Hence it acts on $\sigma_\Lambda^{\overline{I}_{M_\Lambda,1}}$ by the scalar $\omega_{\sigma_1}(\check{\alpha}(t))$.
On the other hand, $\check{\alpha}(t)$ acts on $v$ by the scalar $t^{\langle\nu,\check{\alpha}\rangle} = \nu(\check{\alpha}(t))$.
This gives the lemma.
\end{proof}

\begin{rem}
If we treat the Satake transform in a natural way (see Remark~\ref{rem:choice of uniformizer and Hecke algebra}), Lemma~\ref{lem:Satake parameter of I(lambda)} should be $\Satakepar(I(\Lambda)) = \{(M_{\Pi_1},\omega_{\sigma_1})\}$. (We use a notation of Herzig~\cite[Proposition~4.1]{arXiv:1005.1713}.)
Hence the above lemma should be a consequence of Lemma~\ref{lem:Satake parameter of I(lambda)}
\end{rem}

\begin{prop}\label{prop:irreducibility of I(Lambda)}
For $\Lambda\in\mathcal{P}$, $I(\Lambda)$ is irreducible.
\end{prop}
\begin{proof}
Take $\Lambda = (\Pi_1,\Pi_2,\sigma_1)\in \mathcal{P}$ and put $M_1 = M_{\Pi_1}$ and $M_2 = M_{\Pi_2}$.
Let $\chi$ be an algebra homomorphism $\overline{\kappa}[X_{*,+}]\to\overline{\kappa}$ corresponding to $(M_1,\chi_{\omega_{\sigma_1}})$.
Then $\Satakepar(I(\Lambda)) = \{\chi\}$.
Let $\pi\subset I(\Lambda)$ be a subrepresentation of $I(\Lambda)$.
Take an irreducible $K$-subrepresentation $V$ of $\pi$.
Then $\emptyset\ne \Satakepar(\pi,V)\subset \Satakepar(I(\Lambda)) = \{\chi\}$.
Therefore, we have a nonzero homomorphism $\cInd_K^G(V)\otimes_{\mathcal{H}_G(V)}\chi\to \pi$.

Let $\nu$ be a lowest weight of $V$.
We take $V$ such that the set $\{\alpha\in \Pi\setminus \Pi_{M_\Lambda}\mid \langle \nu,\check{\alpha}\rangle = 0\}$ is minimal.
We claim that this set is empty.
Assume that there exists $\alpha\in\Pi\setminus\Pi_{M_\Lambda}$ such that $\langle\check{\alpha},\nu\rangle = 0$.
Put $\nu' = \nu - (q - 1)\omega_\alpha$ and let $V'$ be the irreducible $K$-representation with lowest weight $\nu'$.
Since $\alpha\not\in\Pi_{M_\Lambda}$, we have $\alpha\not\in\Pi_{\sigma_1}$.
By the definition of $\Pi_{\sigma_1}$, we have:
\begin{itemize}
\item $\langle\check{\alpha},\Pi_{M_1}\rangle \ne 0$ or
\item $\omega_{\sigma_1}(\check{\alpha}(\varpi)) \ne 1$ or $\omega_{\sigma_1}\circ\check{\alpha}|_{\mathcal{O}^\times}$ is not trivial.
\end{itemize}
The above lemma shows that if $\langle\check{\alpha},\Pi_{M_1}\rangle = 0$ then $\omega_{\sigma_1}\circ\check{\alpha}|_{\mathcal{O}^\times}$ is trivial.
Therefore we have that $\langle\check{\alpha},\Pi_{M_1}\rangle \ne 0$ or $\chi_{\omega_{\sigma_1}}(\check{\alpha}) \ne 1$.
Hence we have $\cInd_K^G(V)\otimes_{\mathcal{H}_G(V)}\chi\simeq\cInd_K^G(V')\otimes_{\mathcal{H}_G(V')}\chi$ by Theorem~\ref{thm:changing the weight}.
Therefore, we get a nonzero homomorphism $\cInd_K^G(V')\otimes_{\mathcal{H}_G(V')}\chi\to \pi$.
Namely, $V'$ is an irreducible $K$-subrepresentation of $\pi$.
This contradicts to the minimality of $\{\alpha\in \Pi\setminus \Pi_{M_\Lambda}\mid \langle \check{\alpha},\nu\rangle = 0\}$.

Therefore, we have $\langle \nu,\check{\alpha}\rangle \ne 0$ for $\alpha\in\Pi\setminus \Pi_{M_\Lambda}$.
Put $V_1 = V^{\overline{N_\Lambda}(\kappa)}$.
By \cite[Theorem~3.1]{arXiv:1005.1713}, $\cInd_K^V(G)\otimes_{\mathcal{H}_G(V)}\chi\simeq\Ind_{P_\Lambda}^G(\cInd_{M_\Lambda\cap K}^{M_\Lambda}(V_1)\otimes_{\mathcal{H}_{M_\Lambda}(V_1)}\chi)$.
Therefore, we have a homomorphism $\Ind_{P_\Lambda}^G(\cInd_{M_\Lambda\cap K}^{M_\Lambda}(V_1)\otimes_{\mathcal{H}_{M_\Lambda}(V_1)}\chi)\to\pi\hookrightarrow\Ind_{P_\Lambda}^G\sigma_\Lambda$.
By Lemma~\ref{lem:hom between parabolic induction}, the composition is given by a homomorphism $\cInd_{M_\Lambda\cap K}^{M_\Lambda}(V_1)\otimes_{\mathcal{H}_{M_\Lambda}(V_1)}\chi\to \sigma_\Lambda$.
Since $\sigma_\Lambda$ is irreducible, this homomorphism is surjective.
Therefore, $\cInd_K^G(V)\otimes_{\mathcal{H}_G(V)}\chi\to \Ind_{P_\Lambda}^G(\sigma_\Lambda)$ is surjective.
In particular, $\pi\hookrightarrow \Ind_{P_\Lambda}^G(\sigma_\Lambda)$ is surjective.
Hence $\pi = \Ind_{P_\Lambda}^G(\sigma_\Lambda)$.
\end{proof}

\subsection{Classification theorem}
We will use the following lemma.
\begin{lem}\label{lem:subquotient of parabolic induction from supersingular}
Let $P = MN$ be a parabolic subgroup, $\sigma$ an irreducible admissible representation of $M$ which is supersingular with respect to $(K,T,B)$ and $\omega_\sigma$ the central character of $\sigma$.
Then $\Ind_P^G(\sigma)$ has a filtration whose graded pieces are are $\{I(\Pi_M,\Pi_2,\sigma)\mid \Pi_2\subset\Pi_{\sigma}\}$.
\end{lem}
\begin{proof}
Let $P' = M'N'$ be the standard parabolic subgroup corresponding to $\Pi_M\cup \Pi_{\sigma}$.
Then by Lemma~\ref{lem:G/[M_2,M_2]simeq M_1/L_2}, we can extend $\sigma$ to $M'$ such that $[M_{\Pi_{\sigma}},M_{\Pi_{\sigma}}]$ acts on it trivially.
We have $\Ind_{P\cap M'}^{M'}(\sigma) = (\Ind_{P\cap M'}^{M'}\trivrep_{M})\otimes\sigma$.
So we have $\Ind_P^G(\sigma) = \Ind_{P'}^G((\Ind_{P\cap M'}^{M'}\trivrep_{M'})\otimes\sigma)$.
The definition of the special representations implies that $\Ind_{P\cap M'}^{M'}\trivrep_{M'}$ has a filtration whose graded pieces are $\{\Sp_{Q_2,M'}\}$ where $Q_2$ is a parabolic subgroup of $M'$ which contains $P\cap M'$.
Hence $\Ind_P^G(\sigma)$ has a filtration whose graded pieces are $\{\Ind_{P'}^G(\Sp_{Q_2,M'}\otimes\sigma)\}$.
Let $\Pi_2'\subset\Pi_{M'}$ be a subset corresponding to $Q_2$.
Then we have $\Ind_{P'}^G(\Sp_{Q_2,M'}\otimes\sigma) = I(\Pi_{M},\Pi_2'\setminus\Pi_M,\sigma)$.
\end{proof}
\begin{rem}\label{rem:on composition factors of Ind}
If the derived group of $G$ is simply connected, then $I(\Lambda)$ is irreducible by Proposition~\ref{prop:irreducibility of I(Lambda)}.
Hence the above lemma gives composition factors of $\Ind_P^G(\sigma)$.
In particular, it has a finite length.
The irreducibility of $I(\Lambda)$ will be proved in subsection~\ref{subsec:General case and corollaries}.
Hence the above lemma gives composition factors of $\Ind_P^G(\sigma)$ for any $G$.
\end{rem}

\begin{prop}\label{prop:Main theorem}
Assume that the derived group of $G$ is simply connected.
The correspondence $\Lambda\mapsto I(\Lambda)$ gives a bijection between $\mathcal{P}$ and the set of isomorphism classes of irreducible admissible representations.
\end{prop}
\begin{proof}
First, we prove that the map is surjective  by induction on $\#\Pi$.
Let $\pi$ be an irreducible admissible representation.
Let $\chi$ be an element of $\Satakepar(\pi)$ and assume that it is parameterized by $(M_1,\chi_{M_1})$.
We assume that $M_1$ is minimal.
If $M_1 = G$, then $\pi$ is supersingular.
Therefore, we assume that $M_1\ne G$.
Take an irreducible $K$-representation $V$ such that $\chi\in\Satakepar(\pi,V)$.
Let $\nu$ be a lowest weight of $V$.
We assume that $\Pi_\nu$ is minimal with respect to the condition $\chi\in\Satakepar(\pi,V)$.

Assume that there exists $\alpha\in\Pi_\nu\setminus\Pi_{M_1}$ such that $\langle\Pi_{M_1},\check{\alpha}\rangle \ne 0$ or $\chi_{M_1}(\check{\alpha})\ne 1$.
Set $\nu' = \nu - (q - 1)\omega_\alpha$ and let $V'$ be an irreducible $K$-representation with lowest weight $\nu'$.
Then $\Pi_{\nu'} = \Pi_\nu\setminus\{\alpha\} \subsetneq \Pi_\nu$.
By Theorem~\ref{thm:changing the weight}, we have $\cInd_K^G(V)\otimes_{\mathcal{H}_G(V)}\chi\simeq\cInd_K^G(V')\otimes_{\mathcal{H}_G(V')}\chi$.
Hence $\chi\in\Satakepar(\pi,V')$.
This contradicts to the minimality of $\Pi_\nu$.
Therefore, for all $\alpha\in\Pi_\nu\setminus \Pi_{M_1}$, $\langle\Pi_{M_1},\check{\alpha}\rangle = 0$ and $\chi_{M_1}(\check{\alpha}) = 1$.
From the first condition, $\langle\Pi_\nu\setminus\Pi_{M_1},\check{\Pi}_{M_1}\rangle = 0$.

Let $P = MN$ be a parabolic subgroup corresponding to $\Pi_\nu\cup\Pi_{M_1}$.
First assume that $M\ne G$.
Put $V_1 = V^{\overline{N}(\kappa)}$.
Then we have $\cInd_K^G(V)\otimes_{\mathcal{H}_G(V)}\chi\simeq\Ind_P^G(\cInd_{M\cap K}^M(V_1)\otimes_{\mathcal{H}_M(V_1)}\chi)$~\cite[Theorem~3.1]{arXiv:1005.1713}.
Recall that we have a surjective homomorphism $\cInd_K^G(V)\otimes_{\mathcal{H}_G(V)}\chi\to \pi$.
Hence there exist an irreducible admissible representation $\sigma$ of $M$ and a surjective homomorphism $\Ind_P^G(\sigma)\to \pi$~\cite[Lemma~9.9]{arXiv:1005.1713}.
By inductive hypothesis, $\sigma = I_M(\Lambda')$ for some $\Lambda'\in\mathcal{P}_M$.
Hence there exists a parabolic subgroup $P_0 = M_0N_0$ and an irreducible admissible representation $\sigma_0$ of $M_0$ which is supersingular with respect to $(M_0\cap K,T,M_0\cap B)$ such that $\sigma$ is a subquotient of $\Ind_{P_0\cap M}^M\sigma_0$ by Lemma~\ref{lem:subquotient of parabolic induction from supersingular}.
Hence $\pi$ is a subquotient of $\Ind_{P_0}^G(\sigma_0)$.
By Lemma~\ref{lem:subquotient of parabolic induction from supersingular}, all composition factors of $\Ind_{P_0}^G(\sigma_0)$ are $I(\Lambda)$ for some $\Lambda\in\mathcal{P}$.
Hence $\pi = I(\Lambda)$ for some $\Lambda\in\mathcal{P}$.

Therefore, we may assume that $\Pi_\nu\cup\Pi_{M_1} = \Pi$.
Let $P' = M'N'$ be the standard parabolic subgroup corresponding to $\Pi\setminus\Pi_{M_1}$.
Then for all $\alpha\in\Pi_{M'}$, $\langle\nu,\check{\alpha}\rangle = 0$, $\langle\alpha,\check{\Pi}_{M_1}\rangle = 0$ and $\chi_{M_1}(\check{\alpha}) = 1$.
Set $L' = [M',M']$.
Then the group of coweights $X_{L',*}$ of $L'\cap T$ is $\Z\Pi_{M'}$ which is a subgroup of $X_*\cap\Pi_{M_1}^\perp$.
Put $X_{L',*,+} = X_{*,+}\cap \Z\Pi_{M'}$.
By Lemma~\ref{lem:comparison of satake parameters, restriction to Levi} and Proposition~\ref{prop:comparison of Satake parameters, derived group}, we have $\Satakepar(\pi,V)|_{\overline{\kappa}[X_{L',*,+}]}\subset\Satakepar(\pi|_{M'},V^{\overline{N}'(\kappa)})|_{\overline{\kappa}[X_{L',*,+}]}\subset \Satakepar(\pi|_{L'},V^{\overline{N}'(\kappa)}|_{L'\cap K})$.
Since $\langle\nu,\check{\Pi}_{M'}\rangle = 0$, $V^{\overline{N}'(\kappa)}|_{L'\cap K}$ is trivial.
Therefore, $\chi|_{\overline{\kappa}[X_{L',*,+}]} \in\Satakepar(\pi|_{L'},\trivrep_{L'\cap K})$.
Set $\chi' = \chi|_{\overline{\kappa}[X_{L',*,+}]}$.
We have a non-zero homomorphism $\cInd_{L'\cap K}^{L'}\trivrep_{L'\cap K}\otimes_{\mathcal{H}_{L'}(\trivrep_{L'\cap K})}\chi'\to \pi$.
Since $\chi$ is parameterized by $(M_1,\chi_{M_1})$, $\chi'$ is parameterized by $(L'\cap T,\chi_{M_1}|_{X_{L',*}})$.
Since we have $\chi_{M_1}(\check{\alpha}) = 1$ for all $\alpha\in\Pi_{M'}$, we have $\chi_{M_1}|_{X_{L',*}} = \trivrep_{X_{L',*}}$.
Hence $\chi'$ is parameterized by $(L'\cap T,\trivrep_{X_{L',*}})$.
Therefore, the set of composition factors of $\cInd_{L'\cap K}^{L'}\trivrep_{L'\cap K}\otimes_{\mathcal{H}_{L'}(\trivrep_{L'\cap K})}\chi'$ is $\{\Sp_{Q',L'}\mid \text{$Q'\subset L'$ is a parabolic subgroup}\}$ by Proposition~\ref{prop:changing the weight, by character}.
Hence there exists a parabolic subgroup $P_2 = M_2N_2$ such that $\Pi_{M_1}\subset \Pi_{M_2}$ and $\Sp_{P_2\cap L',L'}\hookrightarrow\pi$.
Let $\sigma_2$ be the special representation $\Sp_{P_2}$.
Then the restriction of $\sigma_2$ to $L'$ is $\Sp_{P_2\cap L',L'}$.
Put $\sigma_1 = \Hom_{L'}(\sigma_2,\pi)$.
This is non-zero.
By Lemma~\ref{lem:general lemma, tensor is irreducible}, $\sigma_1$ is an irreducible representation of $G$ and $\sigma_1\otimes\sigma_2\xrightarrow{\sim}\pi$.

We prove $\sigma_1$ is supersingular as a representation of $M_1$.
Since $L'$ acts on $\sigma_1$ trivially, $\sigma_1$ is regarded as a representation of $G/L'$.
By Lemma~\ref{lem:G/[M_2,M_2]simeq M_1/L_2}, $M_1\to G/L'$ is surjective.
Therefore, $\sigma_1|_{M_1}$ is irreducible.
By inductive hypothesis, $\sigma_1|_{M_1} \simeq I_{M_1}(\Lambda')$.
In particular, $\#\Satakepar(\sigma_1|_{M_1}) = 1$.
Since $\chi$ is parameterized by $(M_1,\chi_{M_1})$, an element of $\Satakepar(\sigma_1|_{M_1})$ is parameterized by $(M_1,\chi_{M_1}')$ for some $\chi_{M_1}'$ by Corollary~\ref{cor:zero-point of Satake parameter of tensor}.
Hence $\sigma_1$ is supersingular.

We prove that the map is injective.
Let $\Lambda' = (\Pi_1',\Pi_2',\sigma_1')$ and assume that $I(\Lambda)\simeq I(\Lambda')$.
Then we have $\Satakepar(I(\Lambda),V) = \Satakepar(I(\Lambda'),V)\ne\emptyset$ for some irreducible representation $V$ of $K$.
By Lemma~\ref{lem:Satake parameter of I(lambda)}, $(M_{\Pi_1},\chi_{\omega_{\sigma_1}}) = (M_{\Pi_1'},\chi_{\omega_{\sigma_1'}})$.
Hence $\Pi_1 = \Pi_1'$.
Let $\nu$ be a lowest weight of $V$.
Then by Lemma~\ref{lem:Stake parameter of I(Lambda), K-part}, for $\alpha\in\Pi$ such that $\langle\Pi_1, \check{\alpha}\rangle = 0$, $\omega_{\sigma_1}\circ \check{\alpha}|_{\mathcal{O}^\times} = \nu\circ\check{\alpha} = \omega_{\sigma_1'}\circ \check{\alpha}|_{\mathcal{O}^\times}$.
On the other hand, we have $\omega_{\sigma_1}\circ\check{\alpha}(\varpi) = \chi_{\omega_{\sigma_1}}(\check{\alpha}) = \omega_{\sigma_1'}\circ\check{\alpha}(\varpi)$.
Hence $\omega_{\sigma_1}\circ\check{\alpha} = \omega_{\sigma_1'}\circ\check{\alpha}$.
Therefore, we have $\Pi_{\sigma_1} = \Pi_{\sigma_1'}$.
Hence $P_\Lambda = P_{\Lambda'}$.

Now we have $\Ind_{P_\Lambda}^G(\sigma_\Lambda)\simeq\Ind_{P_\Lambda}^G(\sigma_{\Lambda'})$.
By Lemma~\ref{lem:hom between parabolic induction}, we have a nonzero homomorphism $\sigma_\Lambda\to\sigma_{\Lambda'}$.
Since $\sigma_\Lambda$ and $\sigma_{\Lambda'}$ are irreducible, $\sigma_\Lambda\simeq\sigma_{\Lambda'}$.
Set $L = [M_{\Pi_{\sigma_1}},M_{\Pi_{\sigma_1}}]$.
As a representation of $L$, $\sigma_\Lambda$ is a direct sum of special representations $\Sp_{Q_2,L}$ where $Q_2$ is a parabolic subgroup of $L$ corresponding to $\Pi_2$.
Hence we have $\Pi_2 = \Pi_2'$.
Therefore, $\sigma_{\Lambda,2} \simeq \sigma_{\Lambda',2}$.
Hence we have $\sigma_1\simeq\Hom_{L}(\sigma_{2,\Lambda},\sigma_\Lambda)\simeq\Hom_L(\sigma_{2,\Lambda'},\sigma_{\Lambda'})\simeq\sigma_1'$.
We get $\Lambda = \Lambda'$.
\end{proof}

\subsection{General case and corollaries}\label{subsec:General case and corollaries}
\begin{thm}\label{thm:Main theorem, general case}
Let $G$ be a connected split reductive algebraic group.
Then $I(\Lambda)$ is irreducible for all $\Lambda\in\mathcal{P}$ and $\Lambda\mapsto I(\Lambda)$ gives a bijection between $\mathcal{P}$ and the set of isomorphism classes of irreducible admissible representations.
\end{thm}
\begin{proof}
Take a $z$-extension $1\to Z\to \widetilde{G}\to G\to 1$ of $G$.
For each parabolic subgroup $P = MN$, let $\widetilde{M}$ be the Levi subgroup of the parabolic subgroup of $\widetilde{G}$ corresponding to $\Pi_M$.
Then $1\to Z\to \widetilde{M}\to M\to 1$ is a $z$-extension of $M$.
For each representation $\pi$ of $G$, let $\widetilde{\pi}$ be the pull-back of $\pi$ to $\widetilde{G}$.
Then we have $I_G(\Pi_1,\Pi_2,\sigma_1)^\sim = I_{\widetilde{G}}(\Pi_1,\Pi_2,\widetilde{\sigma_1})$.
By Proposition~\ref{prop:irreducibility of I(Lambda)}, this is irreducible.
Hence $I_G(\Lambda)$ is irreducible for $\Lambda\in\mathcal{P}$.

We also have that $I_G(\Pi_1,\Pi_2,\sigma_1)\simeq I_G(\Pi_1',\Pi_2',\sigma_1')$ if and only if $I_{\widetilde{G}}(\Pi_1,\Pi_2,\widetilde{\sigma_1})\simeq I_{\widetilde{G}}(\Pi_1',\Pi_2',\widetilde{\sigma_1'})$.
Hence we have $\Pi_1 = \Pi_1'$, $\Pi_2 = \Pi_2'$ and $\widetilde{\sigma_1}\simeq\widetilde{\sigma_1'}$ by Proposition~\ref{prop:Main theorem}.
Hence we have $\sigma_1\simeq \sigma_1'$.

Let $\pi$ be an irreducible admissible representation of $G$.
Then there exists $\Lambda_0 = (\Pi_1,\Pi_2,\sigma_{1,0})\in \mathcal{P}_{\widetilde{G}}$ such that $\widetilde{\pi} = I_{\widetilde{G}}(\Lambda)$.
Since $Z$ is contained in the center of $M_{\Pi_1}$, it acts on $\sigma_{1,0}$ by the scalar.
By the construction of $I_{\widetilde{G}}(\Lambda)$, $Z$ acts on $I_{\widetilde{G}}(\Lambda)\simeq\widetilde{\pi}$ by the same scalar.
It is trivial since $Z$ acts on $\widetilde{\pi}$ trivially.
Hence $Z$ acts on $\sigma_{1,0}$ trivially, namely, $\sigma_{1,0}\simeq\widetilde{\sigma_1}$ for some representation of $G$.
Hence $\pi = I_G(\Pi_1,\Pi_2,\sigma_1)$.
We get the theorem.
\end{proof}
We give corollaries of this theorem.
\begin{cor}
For any irreducible admissible representation $\pi$ of $G$, $\#\Satakepar(\pi) = 1$.
\end{cor}
\begin{proof}
Obvious from Lemma~\ref{lem:Satake parameter of I(lambda)} and Theorem~\ref{thm:Main theorem, general case}.
\end{proof}

\begin{cor}\label{cor:supersinguality and supercuspidality}
Let $\pi$ be an irreducible admissible representation of $G$.
Then the following conditions are equivalent.
\begin{enumerate}
\item The representation $\pi$ is supersingular.
\item The representation $\pi$ is supersingular with respect to $(K,T,B)$.
\item The representation $\pi$ is supercuspidal.
\end{enumerate}
\end{cor}
\begin{proof}
Take $\Lambda = (\Pi_1,\Pi_2,\sigma_1)\in\mathcal{P}$ such that $\pi = I(\Lambda)$.
Then by Lemma~\ref{lem:Satake parameter of I(lambda)}, $\pi$ is supersingular with respect to $(K,T,B)$ if and only if $\Pi_1 = \Pi$.
By Lemma~\ref{lem:subquotient of parabolic induction from supersingular}, $\pi$ is a subquotient of $\Ind_{P_1}^G(\sigma_1)$.
Hence, if $\pi$ is not supersingular with respect to $(K,T,B)$, then $\pi$ is not supercuspidal.

Assume that $\pi$ is a subquotient of $\Ind_{P_0}^G\sigma_0$ for a proper parabolic subgroup $P_0 = M_0N_0$ and an irreducible admissible representation $\sigma_0$.
By Lemma~\ref{lem:subquotient of parabolic induction from supersingular}, we may assume $\sigma_0$ is supersingular with respect to $(K,T,B)$.
By Lemma~\ref{lem:subquotient of parabolic induction from supersingular}, $P_{\Pi_1} = P_0$.
Hence $\pi$ is not supersingular with respect to $(K,T,B)$.

Hence (2) and (3) are equivalent.
Since the property (3) is independent of a choice of $(K,T,B)$, (2) and (1) are equivalent.
\end{proof}

\begin{cor}
Let $P = MN$ be a parabolic subgroup and $\sigma$ a finite length admissible representation of $M$.
Then $\Ind_P^G\sigma$ has a finite length.
\end{cor}
\begin{proof}
We may assume $\sigma$ is irreducible.
This follows from Lemma~\ref{lem:subquotient of parabolic induction from supersingular} and Remark~\ref{rem:on composition factors of Ind},
\end{proof}

\begin{cor}\label{lem:length of principal series}
Let $\nu\colon T\to \overline{\kappa}^\times$ be a character.
Then $\Ind_B^G(\nu)$ has a length $2^C$ where $C = \#\{\alpha\in\Pi\mid \nu\circ\check{\alpha} = \trivrep_{\GL_1}\}$.
In particular, $\Ind_B^G(\nu)$ is irreducible if and only if $\nu\circ\check{\alpha}\ne\trivrep_{\GL_1}$ for all $\alpha\in \Pi$.
\end{cor}
\begin{proof}
Notice that any character of $T$ is supersingular.
Hence this follows from Lemma~\ref{lem:subquotient of parabolic induction from supersingular} and Remark~\ref{rem:on composition factors of Ind}.
\end{proof}

\end{document}